\documentclass[11pt]{tran-l}
\usepackage[ansinew]{inputenc}
\usepackage[centertags]{amsmath}
\usepackage{exscale}
\usepackage{relsize}
\usepackage{amsfonts}
\usepackage{amssymb}
\usepackage{amsthm}
\usepackage{newlfont}

\mathchardef\ordinarycolon\mathcode`\:
  \mathcode`\:=\string"8000
  \begingroup \catcode`\:=\active
    \gdef:{\mathrel{\mathop\ordinarycolon}}
  \endgroup
\def\N{{\Bbb N}}

\def\Q{{\Bbb Q}}
\def\R{{\Bbb R}}
\newtheorem*{mthm}{Main Theorem}
\newtheorem{thm}{Theorem}
\newtheorem{lem}[thm]{Lemma}
\newtheorem{prop}[thm]{Proposition}

\newtheorem{rem}[thm]{Remark}
\newtheorem{Cor}[thm]{Corollary}

\newcommand{\D}{\mathcal{D}}
\newcommand{\B}{\mathcal{B}}
\newcommand{\g}{\mathcal{G}}
\newcommand{\h}{\mathcal{H}}
\newcommand{\aut}{\emph{Aut}}
\newcommand{\Aut}{\mbox{Aut}}
\newcommand{\Out}{\mbox{Out}}

\begin{document}
\title[Flag-transitive Steiner \mbox{$3$-designs}]
{The Classification of Flag-Transitive \\ Steiner
\mbox{$3$-Designs}}

\author{Michael Huber}

\address{Mathematisches Institut der Universit\"{a}t T\"{u}bingen,
Auf der Morgenstelle~10, D-72076~T\"{u}bingen, Germany}

\email{michael.huber@uni-tuebingen.de}



\thanks{The author gratefully acknowledges support by the Deutsche
Forschungsgemeinschaft (DFG)}

\date{April 24, 2003}

\commby{William M. Kantor}


\begin{abstract}
We solve the long-standing open problem of classifying all
\mbox{$3$-$(v,k,1)$ designs} with a flag-transitive group of
automorphisms (cf. A. Delandtsheer, \emph{Geom. Dedicata} \textbf{41}
(1992), p. 147; and \emph{in:} “Handbook of Incidence Geometry”, ed.
by F. Buekenhout, Elsevier Science, Amsterdam, 1995, p. 273; but
presumably dating back to 1965). Our result relies on the
classification of the finite \mbox{$2$-transitive} permutation groups.
\end{abstract}

\maketitle

\section{Introduction}\label{intro}

For positive integers $t \leq k \leq v$ and $\lambda$, we define a
\mbox{\emph{$t$-$(v,k,\lambda)$ design}} to be a finite incidence
structure \mbox{$\D=(X,\B,I)$}, where $X$ denotes a set of
\emph{points}, $\left| X \right| =v$, and $\B$ a set of
\emph{blocks}, $\left| \B \right| =b$, with the properties that
each block $B \in \B$ is incident with $k$ points, and each
\mbox{$t$-subset} of $X$ is incident with $\lambda$ blocks. A
\emph{flag} of $\D$ is an incident point-block pair, that is $x
\in X$ and $B \in \B$ such that $(x,B) \in I$. We consider
automorphisms of $\D$ as pairs of permutations on $X$ and $\B$
which preserve incidence, and call a group \mbox{$G \leq \Aut
(\D)$} of automorphisms of $\D$ \emph{flag-transitive}
(respectively \emph{block-transitive}, \emph{point
$t$-transitive}) if $G$ acts transitively on the flags
(respectively transitively on the blocks, $t$-transitively on the
points) of $\D$. For short, $\D$ is said to be, e.g.,
flag-transitive if $\D$ admits a flag-transitive group of
automorphisms.

We call a \mbox{$t$-$(v,k,1)$ design} a \emph{Steiner
\mbox{$t$-design}} (sometimes this is also known as \emph{Steiner
system}). We note that in this case each block is determined by the
set of points which are incident with it, and thus can be identified
with a $k$-subset of $X$ in a unique way. If furthermore $t<k<v$
holds, then we speak of a \emph{non-trivial} Steiner
\mbox{$t$-design}.

As a consequence of the classification of the finite simple
groups, it has been possible in recent years to characterize
Steiner \mbox{$t$-designs}, mainly for $t=2$, with sufficiently
strong transitivity properties (for an overview,
see~\cite[Sect.\,1,\,2]{Buek1988} and~\cite[Sect.\,2]{Kant1993}).
Probably the most general results have been the classification of
all point \mbox{$2$-transitive} Steiner \mbox{$2$-designs} in 1985
by \mbox{W. M.} Kantor~\cite[Thm.\,1]{Kant1985}, and the almost
complete determination of all flag-transitive Steiner
\mbox{$2$-designs} announced in 1990 by F. Buekenhout, A.
Delandtsheer, J. Doyen, P. B. Kleidman, M. W. Liebeck and J.
Saxl~\cite{Buek1990,Del2001,Lieb1998,Saxl2002} (see
also~\cite[Sect.\,3]{Kant1993} for the incomplete case with a
$1$-dimensional affine group of automorphisms).

Nevertheless, for Steiner \mbox{$3$-designs} such
characterizations have remained challenging open problems. In
particular, the classification of all flag-\linebreak transitive
Steiner \mbox{$3$-designs} is known as “a long-standing and still
open problem” (cf.~\cite[p.\,147]{Del1992}
and~\cite[p.\,273]{Del1995}). Presumably, H.
L\"{u}neburg~\cite{Luene1965} in 1965 has been the first dealing
with part of the problem characterizing flag-transitive Steiner
\mbox{$3$-designs} with block size $k=4$ under the additional
strong assumption that every non-identity element of the group of
automorphisms fixes at most two points. This result has been
generalized recently by the author~\cite{Hu2001}, omitting the
additional assumption. Moreover, in~\cite{Hu2001diss} the author
determined all flag-transitive Steiner \mbox{$3$-designs} with
block size $k \leq 7$.

In this article, we completely classify all flag-transitive
Steiner \mbox{$3$-designs} with arbitrary block size. Our approach
makes use of the classification of the finite
\mbox{$2$-transitive} permutation groups, which in turn relies on
the classification of the finite simple groups. We state our
result:

\medskip

The classification of all non-trivial Steiner \mbox{$3$-designs}
with a flag-transitive group of automorphisms is as follows

\begin{mthm}
Let $\D=(X,\B,I)$ be a non-trivial Steiner \mbox{$3$-design}. Then
\mbox{$G \leq \aut(\D)$} acts flag-transitively on $\D$ if and
only if one of the following occurs:

\medskip

\begin{enumerate}

\item[(1)] $\D$ is isomorphic to the \mbox{$3$-$(2^d,4,1)$} design
whose points and blocks are the points and planes of the affine
space $AG(d,2)$, and one of the following holds:

\medskip

\begin{enumerate}
\item[(i)] $d \geq 3$, and $G \cong AGL(d,2)$,

\medskip

\item[(ii)] $d=3$, and $G \cong AGL(1,8)$ or $A \mathit{\Gamma}
L(1,8)$,

\medskip

\item[(iii)] $d=4$, and $G_0 \cong A_7$,

\medskip

\item[(iv)] $d=5$, and $G \cong A \mathit{\Gamma} L(1,32)$,
\end{enumerate}

\bigskip

\item[(2)] $\D$ is isomorphic to a \mbox{$3$-$(q^e +1,q+1,1)$}
design whose points are the elements of the projective line
\mbox{$GF(q^e) \cup \{\infty\}$} and whose blocks \linebreak are
the images of \mbox{$GF(q) \cup \{\infty\}$} under $PGL(2,q^e)$
(respectively \linebreak $PSL(2,q^e)$, $e$ odd) with a prime power
\mbox{$q \geq 3$}, \mbox{$e \geq 2$}, and the derived design at
any given point is isomorphic to the \mbox{$2$-$(q^e,q,1)$} design
whose points and blocks are the points and lines of $AG(e,q)$, and
\mbox{$PSL(2,q^e) \leq G \leq P \mathit{\Gamma} L (2,q^e)$},

\bigskip

\item[(3)] $\D$ is isomorphic to a \mbox{$3$-$(q+1,4,1)$} design
whose points are the elements of \mbox{$GF(q) \cup \{\infty\}$}
with a prime power \mbox{$q \equiv 7$ $($\emph{mod} $12)$} and
whose blocks are the images of \mbox{$\{0,1,\varepsilon,\infty\}$}
under $PSL(2,q)$, where $\varepsilon$ is a primitive sixth root of
unity in $GF(q)$, and the derived design at any given point is
isomorphic to the Netto triple system $N(q)$, and \mbox{$PSL(2,q)
\leq G \leq P \mathit{\Sigma} L (2,q)$}.

\bigskip

\item[(4)] $\D$ is isomorphic to the Witt \mbox{$3$-$(22,6,1)$}
design, and \mbox{$G \unrhd M_{22}$}.
\end{enumerate}
\end{mthm}

\bigskip

A detailed description of the \emph{Netto triple system} $N(q)$
can be found in~\cite[Sect.\,3]{Deletal1986}.

\bigskip


\section{Definitions and Preliminary Results}\label{Prelim}

If $\D=(X,\B,I)$ is a \mbox{$t$-$(v,k,\lambda)$} design with $t
\geq 2$, and $x \in X$ arbitrary, then the \emph{derived} design
with respect to $x$ is \mbox{$\D_x=(X_x,\B_x, I_x)$}, where $X_x =
X \backslash \{x\}$, \mbox{$\B_x=\{B \in \B: (x,B)\in I\}$} and
$I_x= I \!\!\mid _{X_x \times \; \B_x}$. In this case, $\D$ is
also called an \emph{extension} of $\D_x$. Obviously, $\D_x$ is a
\mbox{$(t-1)$-$(v-1,k-1,\lambda)$} design.

Let $G$ be a permutation group on a non-empty set $X$. For $g \in
G$, let $\mbox{Fix}_X(g)$ denote the set of fixed points of $g$ in
$X$. We call $G$ \emph{semi-regular} if the identity is the only
element that fixes a point of $X$. If additionally $G$ is
transitive, then it is said to be \emph{regular}. If
\mbox{$\{x_1,\ldots,x_m\} \subseteq X$}, let
$G_{\{x_1,\ldots,x_m\}}$ be its setwise stabilizer and
$G_{x_1,\ldots,x_m}$ its pointwise stabilizer (for short, we often
write $G_{x_1 \ldots x_m}$ in the latter case).

For \mbox{$\D=(X,\B,I)$} a Steiner \mbox{$t$-design} with \mbox{$G
\leq \Aut (\D)$}, let $G_B$ denote the setwise stabilizer of a
block $B \in \B$, and for $x \in X$, we define $G_{xB}= G_x \cap
G_B$.

Let $\N$ be the set of positive integers (in this article, $0
\notin \N$). For $ d \in \N$, let $\Phi_d(x)$ denote the $d$-th
cyclotomic polynomial in $\Q[x]$, and for $2 \leq q \in \N$, we
define
\[\Phi_d^*(q)= \frac{1}{f^n} \Phi_d(q),\]
where $f=(d,\Phi_d(q))$ and $f^n$ is the largest power of $f$
dividing $\Phi_d(q)$ if $f \neq 1$, and $n=1$ otherwise
(cf.~\cite[p.\,431]{Her1974}).

Let $m$ and $n$ be integers and $p$ a prime. Then $(m,n)$ is the
greatest common divisor of $m$ and $n$. We write $m \mid n$ if $m$
divides $n$, and $p^m \parallel n$ if $p^m$ divides $n$ but
$p^{m+1}$ does not divide $n$. For $2 \leq q \in \N$, we mean by
$\overline{r} \perp q^n-1$ that $\overline{r}$ divides $q^n-1$ but
not $q^m-1$ for all $1 \leq m < n$.

For any $x \in \R$, let $\lfloor x \rfloor$ (respectively $\lceil
x \rceil$) denote the greatest positive integer which is at most
(respectively the smallest positive integer which is at least)
$x$.

All other notation is standard.

\medskip

\pagebreak

The starting point for our investigation to determine all
flag-transitive Steiner \mbox{$3$-designs} is the following
result.

\begin{prop}\label{flag2trs}
Let $\D=(X,\B,I)$ be a non-trivial Steiner \mbox{$t$-design} with
$t \geq 3$. If \mbox{$G \leq \aut(\D)$} acts flag-transitively on
$\D$, then $G$ also acts point \mbox{$2$-transitively} on $\D$.
\end{prop}

\begin{proof}
Let $x \in X$ arbitrary. As \mbox{$G \leq \Aut(\D)$} acts
flag-transitively on $\D$, obviously $G_x$ acts block-transitively
on the derived Steiner \mbox{$(t-1)$-design} $\D_x$. Since
block-transitivity implies point-transitivity for non-trivial
Steiner \mbox{$t$-designs} with $t \geq 2$ by a theorem of
Block~\cite[Thm.\,2]{Block1965}, $G_x$ also acts
point-transitively on $\D_x$, and the claim follows.
\end{proof}

We note that if $t=2$, then it is elementary that conversely the
point $2$-transitivity of \mbox{$G \leq \Aut(\D)$} implies its
flag-transitivity.

The above proposition allows us to make use of the classification
of all finite \mbox{$2$-transitive} permutation groups, which
itself relies on the classification of all finite simple groups
(cf.~\cite{CSK1976,Gor1982,Her1974,Her1985,Hup1957,Kant1985,Mail1895}).

The list of groups is as follows.

Let $G$ be a finite \mbox{$2$-transitive} permutation group on a
non-empty set $X$. Then $G$ is either of

{\bf (A) Affine Type:} $G$ contains a regular normal subgroup $T$
which is elementary Abelian of order $v=p^d$, where $p$ is a
prime. \mbox{If $a$ divides $d$,} and if we identify $G$ with a
group of affine transformations
\[x \mapsto x^g+u\]
of $V=V(d,p)$, where $g \in G_0$ and $u \in V$, then particularly
one of the following occurs:
\begin{enumerate}

\smallskip

\item[(1)] $G \leq A \mathit{\Gamma} L(1,p^d)$

\smallskip

\item[(2)] $G_0 \unrhd SL(\frac{d}{a},p^a)$, $d \geq 2a$

\smallskip

\item[(3)] $G_0 \unrhd Sp(\frac{2d}{a},p^a)$, $d \geq 2a$

\smallskip

\item[(4)] $G_0 \unrhd G_2(2^a)'$, $d=6a$

\smallskip

\item[(5)] $G_0 \cong A_6$ or $A_7$, $v=2^4$

\smallskip

\item[(6)] $G_0 \unrhd SL(2,3)$ or $SL(2,5)$, $v=p^2$,
$p=5,7,11,19,23,29$ or $59$, or $v=3^4$

\smallskip

\item[(7)] $G_0$ contains a normal extraspecial subgroup $E$ of
order $2^5$, and $G_0/E$ is isomorphic to a subgroup of $S_5$,
$v=3^4$

\smallskip

\item[(8)] $G_0 \cong SL(2,13)$, $v=3^6,$
\end{enumerate}

\smallskip

or

\medskip

{\bf (B) Almost Simple Type:} $G$ contains a simple normal
subgroup $N$, and \mbox{$N \leq G \leq \Aut(N)$}. In particular,
one of the following holds, where $N$ and $v=|X|$ are given as
follows:

\pagebreak

\begin{enumerate}

\smallskip

\item[(1)] $A_v$, $v \geq 5$

\smallskip

\item[(2)] $PSL(d,q)$, $d \geq 2$, $v=\frac{q^d-1}{q-1}$, where
$(d,q) \not= (2,2),(2,3)$

\smallskip

\item[(3)] $PSU(3,q^2)$, $v=q^3+1$, $q>2$

\smallskip

\item[(4)] $Sz(q)$, $v=q^2+1$, $q=2^{2e+1}>2$ \hfill (Suzuki
groups)

\smallskip

\item[(5)] $Re(q)$, $v=q^3+1$, $q=3^{2e+1} > 3$ \hfill (Ree
groups)

\smallskip

\item[(6)] $Sp(2d,2)$, $d \geq 3$, $v = 2^{2d-1} \pm 2^{d-1}$

\smallskip

\item[(7)] $PSL(2,11)$, $v=11$

\smallskip

\item[(8)] $PSL(2,8)$, $v=28$ ($N$ is not \mbox{$2$-transitive)}

\smallskip

\item[(9)] $M_v$, $v=11,12,22,23,24$ \hfill (Mathieu groups)

\smallskip

\item[(10)] $M_{11}$, $v=12$

\smallskip

\item[(11)] $A_7$, $v=15$

\smallskip

\item[(12)] $HS$, $v=176$ \hfill (Higman-Sims group)

\smallskip

\item[(13)] $Co_3$, $v=276$. \hfill (smallest Conway group)
\end{enumerate}

\medskip

For basic properties of the listed groups, we refer, e.g.,
to~\cite{Atlas1985} and~\cite[Ch.\,2,\,5]{KlLi1990}.

We will now indicate some helpful combinatorial tools on which we
rely in the sequel. Let $r$ (respectively $\lambda_2$) denote the
total number of blocks incident with a given point (respectively
pair of distinct points), and let all further parameters be as
defined at the beginning of Section~\ref{intro}.

Obvious is the subsequent fact.

\begin{lem}\label{divprop}
Let $\D=(X,\B,I)$ be a Steiner \mbox{$t$-design}. If $G \leq
\aut(\D)$ acts flag-transitively on $\D$, then, for any $x \in X$,
the division property
\[r  \bigm|  \left| G_x \right|\] holds.
\end{lem}

Elementary counting arguments give the following standard assertions.

\begin{lem} \label{Comb_t=3}
If $\D=(X,\B,I)$ is a \mbox{$t$-$(v,k,\lambda)$} design, then the
following holds:
\begin{enumerate}

\smallskip

\item[(a)] $bk = vr.$

\smallskip

\item[(b)] $\displaystyle{{v \choose t} \lambda = b {k \choose
t}.}$

\smallskip

\item[(c)] $r(k-1)=\lambda_2(v-1)$ for $t \geq 2$, where
$\displaystyle{\lambda_2=\lambda \frac{{v-2 \choose t-2}}{{k-2
\choose t-2}}.}$

\smallskip

\item[(d)] In particular, if $t=3$, then $(k-2) \lambda_2 = v-2.$
\end{enumerate}
\end{lem}

\medskip
\pagebreak

For non-trivial Steiner \mbox{$t$-designs} lower bounds for $v$ in
terms of $k$ and $t$ can be indicated.

\begin{prop}{\em (Cameron~\cite{Cam1976}).}\label{Cam}
Let $\D=(X,\B,I)$ be a non-trivial Steiner \mbox{$t$-design}. Then
the following holds:
\begin{enumerate}

\smallskip

\item[(a)] $v\geq (t+1)(k-t+1).$

\smallskip

\item[(b)] $v-t+1 \geq (k-t+2)(k-t+1)$ for $t>2$. If equality
holds, then
\smallskip
$(t,k,v)=(3,4,8),(3,6,22),(3,12,112),(4,7,23)$, or $(5,8,24)$.
\end{enumerate}
\end{prop}

We note that (a) is stronger for $k<2(t-1)$, while (b) for
$k>2(t-1)$. For $k=2(t-1)$ both assert that $v \geq t^2-1$.

As we are in particular interested in the case when $t=3$, we
deduce from (b) the following upper bound for the positive integer
$k$.

\begin{Cor}\label{Cameron_t=3}
Let $\D=(X,\B,I)$ be a non-trivial Steiner \mbox{$3$-design}. Then
the block size $k$ can be estimated by
\[k \leq \bigl\lfloor \sqrt{v} + \textstyle{\frac{3}{2}} \bigr\rfloor.\]
\end{Cor}

\smallskip

\begin{rem} \label{equa_t=3}
\emph{If \mbox{$G \leq \Aut(\D)$} acts flag-transitively on any
Steiner \mbox{$3$-design} $\D$, then applying
Proposition~\ref{flag2trs} and Lemma~\ref{Comb_t=3}~(b) yields the
equation
\[b=\frac{v(v-1)(v-2)}{k(k-1)(k-2)}=\frac{v(v-1) \left|G_{xy}\right|}
{\left| G_B \right|},\] where $x$ and $y$ are two distinct points
in $X$ and $B$ is a block in $\B$, and thus \[v-2 = (k-1)(k-2)
\frac{\left|G_{xy} \right|}{\left|G_{xB}\right|} \quad \mbox{if}
\;\, x \in B.\]}
\end{rem}

\bigskip


\section{Cases with a Group of Automorphisms of Affine Type}\label{affine typ}

In the following, we begin with the proof of the Main Theorem.
Using the notation as before, let $\D=(X,\B,I)$ be a non-trivial
Steiner \mbox{$3$-design} with \mbox{$G \leq \Aut(\D)$} acting
flag-transitively on $\D$. Let us recall that in view of
Proposition~\ref{flag2trs}, we can restrict ourselves to the
inspection of the finite \mbox{$2$-transitive} permutation groups
listed in Section~\ref{Prelim}. Before we consider in this section
successively those cases where $G$ is of affine type, we prove
some lemmas which will be required for Case (1).


\begin{lem}\label{AGL1}
Let $q=p^d$ with $p \neq 2$ a prime. Furthermore, let $2^m
\parallel p-1$, $2^{\overline{m}}\parallel p+1$ and $2^n \parallel
d$ for some integers $m$, $\overline{m}$ and $n$. Then $2^{m+n}
\parallel q-1$, unless \mbox{$p \equiv 3$ $($\emph{mod} $4)$} and
$d \equiv 0$ $($\emph{mod} $2)$, in which case $2^{\overline{m}+n}
\parallel q-1.$
\end{lem}
\begin{proof}
This follows from~\cite[Lemma\,3.2]{Her1974} using induction over
$n$.
\end{proof}

\pagebreak

Maintaining the same parameters, we obtain

\begin{lem}\label{AGL2}
Let $G \leq A \mathit{\Gamma} L(1,q)$ be a \mbox{$2$-transitive}
permutation group, where $q=p^d$ with $p \neq 2$ a prime, and $P$
a Sylow $2$-subgroup of $G$. Then we have \mbox{$\left| P \cap
AGL(1,q) \right| \geq 2^m$}. Moreover, if \mbox{$p \equiv 3$
$($\emph{mod} $4)$} and $d \equiv 0$ $($\emph{mod} $2)$, then
\mbox{$\left| P \cap AGL(1,q) \right| \geq 2^{\overline{m}}$}.
\end{lem}
\begin{proof}

Clearly,
\[P / P \cap AGL(1,q) \cong P \cdot AGL(1,q) / AGL(1,q) \leq
A\mathit{\Gamma} L(1,q) / AGL(1,q).\] Thus, we obtain \[\left| P
\right| \bigm|  \left| P \cap AGL(1,q) \right| \cdot d.\] As
$q(q-1) \bigm| \left| G \right|$ by the $2$-transitivity of $G$,
Lemma~\ref{AGL1} yields
\[2^{m+n} \bigm| \left| P \right| \bigm| \left| P \cap AGL(1,q) \right| \cdot 2^n,\]
and therefore
\[2^m \bigm| \left|\, P \cap AGL(1,q) \right|.\]
If $p \equiv 3$ (mod $4$) and $d \equiv 0$ (mod $2$), then we have
$2^{\overline{m}+n} \bigm|  q-1$, and hence $2^{\overline{m}} \bigm| \left| P \cap
AGL(1,q) \right|.$
\end{proof}

\begin{lem}\label{AGL3}
Let $G \leq A \mathit{\Gamma} L(1,q)$ be a \mbox{$2$-transitive}
permutation group, where $q=p^d$ with $p \neq 2$ a prime. Then $G$
contains an involution which fixes exactly one point.
\end{lem}
\begin{proof}
Clearly, $AGL(1,q)_0$ is isomorphic to $GL(1,q)$, and hence
cyclic. It has index $q$, which is odd, and contains therefore a
Sylow $2$-subgroup of $AGL(1,q)$. Thus, each involution in
$AGL(1,q)$ has exactly one fixed point, and the claim follows by
applying Lemma~\ref{AGL2}.
\end{proof}

\medskip

We shall now turn to the examination of those cases where \mbox{$G
\leq \Aut(\D)$} is of affine type.

\bigskip
\emph{Case} (1): $G \leq A \mathit{\Gamma} L(1,v)$, $v=p^d$.
\medskip

First, we will show by contradiction that $v$ is a power of $2$.
Indeed, we suppose that $p \neq 2$. Let $T$ denote the translation
subgroup of $G$. By Lemma~\ref{AGL3}, we know that $G$ contains an
involution $\tau$ which has exactly one fixed point $x \in X$.
Then, for distinct $x,y \in X$, the \mbox{$3$-subset}
$S=\{x,y,y^\tau\}$ is invariant under $\tau$. But, $S$ is incident
with a unique block $B \in \B$ by the definition of Steiner
\mbox{$3$-designs}, hence $\tau \in G_B$. Since $G$ is
flag-transitive, $G_B$ acts transitively on the points of $B$.
Therefore, for each point $x \in B$, there exists an involution
$\tau_x$ having only $x$ as fixed point. Hence
\[U := \text{\footnotesize{$\langle$}}
{\tau_x}^{G_B}\text{\footnotesize{$\rangle$}} \leq
\text{\footnotesize{$\langle$}} {\tau_x}^{A \mathit{\Gamma}
L(1,v)} \text{\footnotesize{$\rangle$}}
=\text{\footnotesize{$\langle$}} \tau_x
\text{\footnotesize{$\rangle$}} \cdot T,\] whereas for the latter
we use that $\tau_x$ induces on $T$ the inverse map
\mbox{$\alpha:x \mapsto x^{-1}$} because any
involutory automorphism of $T$ which has no fixed point distinct
from $1$ must be equal to $\alpha$. Therefore, we have
\mbox{$\tau_x \in AGL(1,v) \trianglelefteq A \mathit{\Gamma}
L(1,v)$}. Then, by Dedekind's law,
\[U= \text{\footnotesize{$\langle$}} \tau_x
\text{\footnotesize{$\rangle$}} \cdot (U \cap T).\]
But, as $U$ acts transitively on the points of $B$ and clearly
$\text{\footnotesize{$\langle$}} \tau_x \text{\footnotesize{$\rangle$}} \cap (U \cap T) = 1$,
it follows from the orbit-stabilizer property that
$U \cap T$ acts also transitively on the points of $B$. Thus, $B$ is a
point-orbit under $U \cap T$ and therefore a subspace of
$AG(d,p)$. Since $G$ is block-transitive, we conclude that all blocks
must be affine subspaces.\\
Let $\g$ be a line in $AG(d,p)$ with distinct points $x,y \in \g$.
Let $B$ and $\overline{B}$ be two distinct blocks containing
$\{x,y\}$. As $p \neq 2$ and since affine subspaces
contain with any two distinct points also the line connecting
them, it follows that \mbox{$\g \subseteq B \cap \overline{B}$}
with $\left| \g \right|> 2$, a contradiction. Thus, we have shown
that $v=2^d$.

In the following, we will prove that if the block size $k$ is a
power of $2$, then only $k=4$ can occur. Therefore, we can use the
classification of all flag-transitive Steiner quadruple
systems~\cite{Hu2001}, which gives the designs described in part
(1) of the Main Theorem with the assertions (ii) and (iv). To
exclude trivial Steiner \mbox{$3$-designs}, let $k=2^a$, $1 <a <
d$. As $d=3$ yields $k=4$, we may assume that $d>3$. From
Remark~\ref{equa_t=3}, it follows that
\begin{equation}\label{eq-1}
v-2 \bigm| d (k-1)(k-2).
\end{equation}
Combining this with~\cite[Thm.\,3.3\,(a)]{Her1974} gives
\begin{equation} \label{eq-1.1}
\Phi_{d-1}^*(2)  \bigm|  2^{d-1}-1  \bigm|  d (2^a-1)(2^{a-1}-1).
\end{equation}
Clearly, $a < d-1$ (otherwise, $k=2^{d-1}$, a contradiction to
Corollary~\ref{Cameron_t=3}.) If $\Phi_{d-1}^*(2)=1$, then,
by~\cite[Thm.\,3.5\,]{Her1974}, there exists no non-trivial
\mbox{$2$-primitive} prime divisor of $2^{d-1}-1$, and hence $d=7$
in view of Zsigmondy's theorem (see~\cite[p.\,283]{Zsig1892}). By
using~(\ref{eq-1}), Lemma~\ref{Comb_t=3}~(d) and
Corollary~\ref{Cameron_t=3}, we can easily check the very small
number of possibilities for $k$. It turns out that only $k=4$ can
occur. Thus, we may assume that there exists a prime divisor
$\overline{r}$ of $\Phi_{d-1}^*(2)$. Then $\overline{r} \mid d$
by~\cite[Thm.\,3.5\,(vi)]{Her1974}. As \mbox{$\overline{r} \equiv
1$ (mod $(d-1)$)} (which follows
from~\cite[Thm.\,3.5\,(ii)]{Her1974}), we conclude that
$\overline{r}=d$. If there exists a further prime divisor
$\hat{r}$ of $\Phi_{d-1}^*(2)$ with $\hat{r} \neq \overline{r}$,
then again $\hat{r} \mid d$ and $\hat{r} =d$ by the same
arguments. Thus $\hat{r}= \overline{r}$, a contradiction. Hence,
we have
\[\Phi_{d-1}^*(2)= \overline{r}^n\]
for some $n \in \N$. But then, by dividing~(\ref{eq-1.1}) by
$\overline{r}$ and using~\cite[Thm.\,3.5\,(vi)]{Her1974} again, we
obtain
\[\frac{\Phi_{d-1}^*(2)}{\overline{r}} \leq 1.\]
Therefore, $\Phi_{d-1}^*(2) \leq \overline{r}=d$. As
$\Phi_{d-1}^*(2)=1$ has already been considered, we may suppose
that $\Phi_{d-1}^*(2)=d$. Now~\cite[Thm.\,3.9\,(b)]{Her1974}
yields $d \leq 19$. The small number of cases can easily be
checked by hand as above. Again, it turns out that only $k=4$ can
occur.

Let us suppose now that $k$ is no power of $2$. We distinguish two
cases according as some non-trivial translation preserves a block
$B \in \B$ or not. Let $T_B \neq 1$. Then $B$ is a disjoint union
of affine subspaces $X_i$ of $AG(d,2)$, $i \geq 1$ (namely the
point-orbits $X_i$ of $T_B$ contained in $B$). As $k$ is no power
of $2$, we may assume that $i \geq 2$. Let $x_i \in X_i$. Then the
translation $t$ mapping $x_1$ onto $x_i$ maps $B$ onto some other
block $B_i$ (because $t \notin T_B$). Since $X_i \subseteq B \cap
B_i$ and $\left| X_i \right| \geq p = 2$, it follows from the
definition of Steiner \mbox{$3$-designs} that $\left| X_i \right|
= 2$ for each $i$. Therefore, $\left| G_B \cap T \right| = \left|
T_B \right| = 2$. Without restriction, we may assume that
$T_B=\text{\footnotesize{$\langle$}} x \mapsto x+1
\text{\footnotesize{$\rangle$}}$. Thus
\[G_B \leq \mathcal{C}_{A \mathit{\Gamma} L(1,v)}(T_B)= T \cdot
\text{\footnotesize{$\langle$}}\alpha
\text{\footnotesize{$\rangle$}},\] where $\mathcal{C}_{A
\mathit{\Gamma} L(1,v)}(T_B)$ denotes the centralizer of $T_B$ in
$A \mathit{\Gamma} L(1,v)$ and $\alpha$ the Frobenius automorphism
$GF(v) \longrightarrow GF(v)$, \mbox{$x \mapsto x^2$}. Hence
\[G_B / T_B \cong G_B \cdot T / T\] is isomorphic to a subgroup of
\[\mathcal{C}_{A \mathit{\Gamma} L(1,v)}(T_B) / T \cong
\text{\footnotesize{$\langle$}} \alpha
\text{\footnotesize{$\rangle$}}.\] Because of the transitivity of
$G_B$ on the points of $B$, we conclude that \linebreak $k \bigm|
\left|G_B \right| \bigm|  2d.$ Therefore, $v-2 < 4d^3$
by~(\ref{eq-1}), and the small number of possibilities for $k$ can
easily be eliminated by hand using~(\ref{eq-1}) and
Lemma~\ref{Comb_t=3}~(d).\\
Now, let $T_B=1$. We first show that $G_B \leq G_y$ for some $y
\notin B$. Let $G^* = G_B \cap AGL(1,v)$. Then $G^*$ is conjugated
to a subgroup of $G_0$ by Hall's theorem. If $G^*=1$, then $G_B$
is isomorphic to a subgroup of $\text{\footnotesize{$\langle$}}
\alpha \text{\footnotesize{$\rangle$}}$, hence cyclic and
$\left|G_B \right| \bigm| d$. As $G_B$ acts transitively on the
points of $B$, we obtain $k \mid d$, and thus $v-2 < d^3$
by~(\ref{eq-1}). The very few possibilities for $k$ can easily be
ruled out by hand as before. Therefore, $G^* \neq 1$. By
construction, $G^*$ has only the point $0$ as fixed point. Since
$G^* \unlhd G_B$, obviously $G_B$ fixes the set of fixed points of
$G^*$, i.e. the point $0$. Hence $G_B \leq G_0$, and $0 \notin B$
by the flag-transitivity of $G$.\\
As $G$ is point \mbox{$2$-transitive}, we have $\left| G \right|
=v(v-1)a$ with $a \mid d$. Then Remark~\ref{equa_t=3} yields
\begin{equation} \label{eq-2}
v-2 = (k-1)(k-2) \frac{a}{\left|G_{xB}\right|} \quad \mbox{if}
\;\, x \in B.
\end{equation}
As $G_B$ fixes some $y \notin B$, it
follows that $\left| G_{xB} \right|\bigm| \left| G_{xy} \right| = a$.\\
If $G_{0x}$ fixes three or more distinct points, then $G_{0x}$
would fix some block $\overline{B} \in \B$. Thus, we have $a
\bigm| \left| G_{xB}\right|$, and therefore $v-2= (k-1)(k-2)$.
However, as $d>3$, it follows from Proposition~\ref{Cam}~(b) that
$v-2
> (k-1)(k-2)$, a contradiction. Hence, $G_{0x}$ fixes only $0$ and
$x$. Then $G_{0x}$ must contain a field automorphism of order
$d$, and we conclude that $G= A \mathit{\Gamma} L(1,2^d)$.\\
Let $p$ be a prime divisor of $d$, say $d=ps$. Then $(G_{0x})^p$
fixes at least three distinct points, and hence we have $s \bigm|
\left| G_{xB} \right|$. If there exists a further prime divisor
$\overline{p}$ of $d$ with $\overline{p} \neq p$, then the
quotients $d/\overline{p}$ and $d / p$ both divide the order of
$G_{xB}$ by the flag-transitivity of $G$. Therefore, we obtain $d
\bigm|
\left| G_{xB} \right|$, which gives the contradiction $a=d$ as above.\\
Thus, we have $d =p^n$ for some $n \in \N$, and therefore $p^{n-1}
= s \bigm| \left| G_{xB} \right|.$ Now, it follows that $\left|
G_{xB} \right|=p^{n-1}$, and hence \mbox{$\left| G_B \right| = k
p^{n-1} \bigm| (v-1)p^n.$} This shows that \mbox{$k \mid (v-1)p$}.
If we set $c=(k,p)$, then $c=1$ or $p$, and we obtain
$\frac{k}{c} \bigm| v-1$. Comparing this with
equation~(\ref{eq-2}) yields
\[v-2 = (k-1)(k-2) \frac{p^n}{p^{n-1}},\]
and hence \[-1 \equiv 2p \; \Big(\mbox{mod} \; \frac{k}{c}\Big).\]
Therefore, we have
\[\frac{k}{c} \leq 2p+1,
\] and finally
\[ 2^{p^n}-2 = v-2= (k-1)(k-2)p \leq (2p^2+p-1)(2p^2+p-2)p.\]
This leaves only a small number of cases to check. As \mbox{$k
\bigm| (2^{p^n}-1)p$}, and $k \geq
\Big\lceil\sqrt{\frac{2^{p^n}-2}{p^n}}+\frac{3}{2}\Big\rceil$
by~(\ref{eq-1}), these can again easily be eliminated by hand
using Lemma~\ref{Comb_t=3}~(c) and (d), and
Corollary~\ref{Cameron_t=3}.

\bigskip
\emph{Case} (2): $G_0 \unrhd SL(\frac{d}{a},p^a)$, $d \geq 2a$.
\medskip

In the following, let $e_i$ denote the $i$-th unit vector of the
vector space $V=V(\frac{d}{a},p^a)$, and
$\text{\footnotesize{$\langle$}} e_i
\text{\footnotesize{$\rangle$}}$ the $1$-dimensional vector
subspace spanned by $e_i$. We will show that only the
flag-transitive designs described in part (1) of the Main Theorem
with $d \geq 3$ and \mbox{$G \cong AGL(d,2)$} can occur.

First, let $p^a \neq 2$. For $d=2a$, let
$U=U(\text{\footnotesize{$\langle$}} e_1
\text{\footnotesize{$\rangle$}}) \leq G_0$ denote the subgroup of
all transvections with axis $\text{\footnotesize{$\langle$}} e_1
\text{\footnotesize{$\rangle$}}$. Then $U$ consists of all
elements of the form
\[\begin{pmatrix}
  1 & 0 \\
  c & 1
\end{pmatrix},\; c \in GF(p^a) \; \mbox{arbitrary}.\]
Clearly, $U$ fixes as points only the elements of
$\text{\footnotesize{$\langle$}} e_1
\text{\footnotesize{$\rangle$}}$. Hence, $G_0$ has point-orbits of
length at least $p^a$ outside $\text{\footnotesize{$\langle$}} e_1
\text{\footnotesize{$\rangle$}}$. Now, let \mbox{$x \in
\text{\footnotesize{$\langle$}} e_1
\text{\footnotesize{$\rangle$}}$} be distinct from $0$ and $e_1$.
Obviously, $U$ fixes the unique block $B \in \B$ which is incident
with the \mbox{$3$-subset} $\{0,e_1,x\}$. Thus, if $B$ contains at
least one point outside $\text{\footnotesize{$\langle$}} e_1
\text{\footnotesize{$\rangle$}}$, then we would obtain $k \geq p^a
+3$. But, according to Corollary~\ref{Cameron_t=3}, we have $k
\leq p^a+1$, a contradiction. Therefore, $B$ is contained
completely in $\text{\footnotesize{$\langle$}} e_1
\text{\footnotesize{$\rangle$}}$. Hence, as $G$ is
flag-transitive, we may conclude that each block lies in an affine
line. But, by the definition of Steiner \mbox{$3$-designs}, any
three distinct non-collinear points must also be incident with a
unique block, a contradiction.\\
For $d \geq 3a$, we consider $(\frac{d}{a} \times
\frac{d}{a})$-matrices of the form
\[A_i = \begin{pmatrix}
  1 & 0 & \hspace{-0.2cm}0 & \hspace{-0.3cm} \cdots & 0 \\
  \vspace{0.1cm}x_1 &  &  &  &  \\
  0 &  & \hspace{0.2cm}\text{\raisebox{-1ex}{\large{{$B_i$}}}} &  &  \\
  \vspace{0.25cm}\vdots &  &  & \\
  0 &  &  &  &
\end{pmatrix}, \, 1 \leq i \leq \textstyle{\frac{d}{a}-1}, \; x_1 \in GF(p^a)\; \mbox{arbitrary},\]
where \[B_1=\begin{pmatrix}
  x_2 & x_3 & x_4 & \cdots \; & x_{\frac{d}{a}}  \\
  0 & \hspace{0.1cm}x_2^{-1} & &    &  \\
  0 &  & \hspace{-0.15cm}1 & \vspace{0.1cm}\hspace{0.5cm}\text{\huge{$*$}} &  \\
  \vspace{0.2cm}\vdots &   \hspace{0.6cm}\text{\Large{$0$}} &   & \hspace{-0.2cm}\ddots \\
  \vspace{0.1cm}0 &  &  &  & 1
\end{pmatrix}, \, x_2 \neq 0,\] \[B_2=\begin{pmatrix}
  0 & x_3 & x_4 & x_5 & \cdots& x_{\frac{d}{a}} \\
  x_3^{-1} & 0 &  &  &  &  \\
  0 &  & \hspace{-0.2cm}-1 &  & \hspace{0.1cm}\text{\huge{$*$}} &  \\
  0 &  &  & 1 &  &  \\
  \vspace{0.1cm}\vdots &  & \hspace{-0.5cm}\text{\Large{$0$}}   & &\ddots &  \\
  \vspace{0.1cm}0 &  &  &  &  & 1
\end{pmatrix}, \, x_3 \neq 0,\]
\[\mbox{and} \;\, B_i=\begin{pmatrix}
  0 & 0 & \cdots & 0 & \hspace{0.2cm}0 & x_{i+1} & x_{i+2}   & \cdots &  x_{\frac{d}{a}}\\
  0 & 1 &  &  &  &  & & &  \\
  \vspace{0.1cm}\vdots &  & \ddots &  &  & & & &  \\
  \vspace{0.1cm}0 &  &  & 1 &  & & \hspace{-0.2cm}\text{\huge{$*$}}& &   \\
  0 &  &  &  & -1 &  & &   \\
  x_{i+1}^{-1} &  &  &  &  & 0 &  &  &   \\
  0 &  &  &  &  &  & \hspace{-0.2cm}1 &  &  \\
  \vspace{0.15cm}\vdots &  &  & \hspace{-0.4cm}\text{\Large{$0$}} & &    & & \hspace{-0.25cm}\ddots & \\
  0 &  &  &  &  &  &  &  &  1
\end{pmatrix}, \, x_{i+1} \neq 0, \; 3 \leq i \leq \textstyle{\frac{d}{a}}-1.\]
Obviously, $B_i \in SL(\frac{d}{a}-1,p^a)$ for $1 \leq i \leq
\frac{d}{a}-1$, and hence $A_i \in SL(\frac{d}{a},p^a)_{e_1}$ by
Laplace's expansion theorem. By multiplying $e_2$ with the
matrices $A_i$ $(1 \leq i \leq \frac{d}{a}-1)$, we obtain as
images exactly the vectors of \mbox{$V \setminus
\text{\footnotesize{$\langle$}} e_1
\text{\footnotesize{$\rangle$}}$}. Thus
$SL(\frac{d}{a},p^a)_{e_1}$, and hence also $G_{0,e_1}$, acts
point-transitively outside $\text{\footnotesize{$\langle$}} e_1
\text{\footnotesize{$\rangle$}}$. Again, let $x \in
\text{\footnotesize{$\langle$}} e_1
\text{\footnotesize{$\rangle$}}$ be distinct from $0$ and $e_1$.
If the unique block $B \in \B$ which is incident with the
\mbox{$3$-subset} $\{0,e_1,x\}$ contains some point outside
$\text{\footnotesize{$\langle$}} e_1
\text{\footnotesize{$\rangle$}}$, then it would already contain
all points outside, thus at least $p^d-p^a+3$ many, which
obviously contradicts Corollary~\ref{Cameron_t=3}. Therefore, $B$
lies completely in $\text{\footnotesize{$\langle$}} e_1
\text{\footnotesize{$\rangle$}}$, and by the same argument as
above, we obtain that here \mbox{$G \leq \Aut(\D)$} cannot act
flag-transitively on any non-trivial Steiner \mbox{$3$-design}
$\D$.

Now, let $p^a=2$. To obtain non-trivial Steiner
\mbox{$3$-designs}, let $v=2^d > 4$. For $v=8$, necessarily $k=4$
must hold in view of Lemma~\ref{Comb_t=3}~(c). For $v>8$, we will
show that also only Steiner quadruple systems can occur. Thus,
applying~\cite{Hu2001} yields the claim. We remark that clearly
any three distinct points are non-collinear in $AG(d,2)$ and hence
define an affine plane. Let $\mathcal{E}=
\text{\footnotesize{$\langle$}} e_1,e_2
\text{\footnotesize{$\rangle$}}$ denote the \mbox{$2$-dimensional}
vector subspace spanned by $e_1$ and $e_2$. We consider $(d \times
d)$-matrices of the form
\[A_i=\begin{pmatrix}
  1 & 0 & \hspace{0.15cm}0 & \hspace{0.05cm}\cdots & 0  \\
  0 & 1 & \hspace{0.15cm}0 & \hspace{0.05cm}\cdots  & 0 \\
  x_1 & x_2 &  &  &   \\
  \vspace{-0.15cm}0 & 0 &  & \hspace{-0.4cm}\text{\large{$B_i$}} &  \\
  \vdots & \vdots &  &  &   \\
  0 & 0 &  &  &
\end{pmatrix}, \, 1 \leq i \leq d-2;\; x_1,x_2 \in GF(2)\; \mbox{arbitrary}\]
with
\[B_1=\begin{pmatrix}
  x_3 & x_4 & x_5 & \cdots \; & x_d  \\
  0 & \hspace{0.1cm}x_3^{-1} & &    &  \\
  0 &  & \hspace{-0.15cm}1 & \vspace{0.1cm}\hspace{0.5cm}\text{\huge{$*$}} &  \\
  \vspace{0.2cm}\vdots &   \hspace{0.6cm}\text{\Large{$0$}} &   & \hspace{-0.2cm}\ddots \\
  \vspace{0.1cm}0 &  &  &  & 1
\end{pmatrix}, \, x_3 \neq 0,\] \[B_2=\begin{pmatrix}
  0 & x_4 & x_5 & x_6 & \cdots& x_d \\
  x_4^{-1} & 0 &  &  &  &  \\
  0 &  & \hspace{-0.2cm}-1 &  & \hspace{0.1cm}\text{\huge{$*$}} &  \\
  0 &  &  & 1 &  &  \\
  \vspace{0.1cm}\vdots &  & \hspace{-0.5cm}\text{\Large{$0$}}   & &\ddots &  \\
  \vspace{0.1cm}0 &  &  &  &  & 1
\end{pmatrix}, \, x_4 \neq 0,\]
\[\mbox{and} \;\, B_i=\begin{pmatrix}
  0 & 0 & \cdots & 0 & \hspace{0.2cm}0 & x_{i+2} & x_{i+3}   & \cdots &  x_d\\
  0 & 1 &  &  &  &  & & &  \\
  \vspace{0.1cm}\vdots &  & \ddots &  &  & & & &  \\
  \vspace{0.1cm}0 &  &  & 1 &  & & \hspace{-0.2cm}\text{\huge{$*$}}& &   \\
  0 &  &  &  & -1 &  & &   \\
  x_{i+2}^{-1} &  &  &  &  & 0 &  &  &   \\
  0 &  &  &  &  &  & \hspace{-0.2cm}1 &  &  \\
  \vspace{0.15cm}\vdots &  &  & \hspace{-0.4cm}\text{\Large{$0$}} & &    & & \hspace{-0.25cm}\ddots & \\
  0 &  &  &  &  &  &  &  &  1
\end{pmatrix}, \, x_{i+2} \neq 0, \; 3 \leq i \leq d-2.\]
Analogously as above, $B_i \in SL(d-2,2)$ for $1 \leq i \leq d-2$
and \mbox{$A_i \in SL(d,2)_\mathcal{E}$}. By multiplying $e_3$
with the matrices $A_i$ $(1 \leq i \leq d-2)$, we obtain as images
exactly the vectors of \mbox{$V \setminus \mathcal{E}$}. Hence
$SL(d,2)_\mathcal{E}$, and therefore also $G_{0,\mathcal{E}}$,
acts point-transitively on \mbox{$V \setminus \mathcal{E}$}. If
the unique block $B \in \B$ which is incident with the
\mbox{$3$-subset} $\{0,e_1,e_2\}$ contains some point outside
$\mathcal{E}$, then it would already contain all points of
\mbox{$V \setminus \mathcal{E}$}. But then, we would have \mbox{$k
\geq 2^d-4+3=2^d-1$}, a contradiction to
Corollary~\ref{Cameron_t=3}. Hence, $B$ lies completely in
$\mathcal{E}$, and by the flag-transitivity of $G$, it follows
that each block must be contained in an affine plane. Thus $k \leq
4$, and finally $k=4$ as we exclude trivial Steiner
\mbox{$3$-designs}.

\bigskip
\emph{Case} (3): $G_0 \unrhd Sp(\frac{2d}{a},p^a)$, $d \geq 2a$.
\medskip

We will prove by contradiction that \mbox{$G \leq \Aut(\D)$}
cannot act flag-\linebreak transitively on any non-trivial Steiner
\mbox{$3$-design} $\D$. First, let $p^a \neq 2.$ The permutation
group $PSp(\frac{2d}{a},p^a)$ on the points of the associated
projective space is a \mbox{rank $3$} group, and the orbits of the
one-point stabilizer are known
(e.g.~\cite[Ch.\,II,\,Thm.\,9.15\,(b)]{HupI1967}). Thus, $G_0
\unrhd Sp(\frac{2d}{a},p^a)$ has exactly two orbits on \mbox{$V
\setminus \text{\footnotesize{$\langle$}} x
\text{\footnotesize{$\rangle$}}$} $(0 \neq x \in V)$ of length at
least
\[\frac{p^a(p^{2d-2a}-1)}{p^a-1}=\sum_{i=1}^{\frac{2d}{a}-2}p^{ia}
> p^d.\] Let $y \in \text{\footnotesize{$\langle$}} x
\text{\footnotesize{$\rangle$}}$ be distinct from $0$ and $x$. If
the unique block which is incident with the \mbox{$3$-subset}
$\{0,x,y\}$ contains at least one point of \mbox{$V \setminus
\text{\footnotesize{$\langle$}} x
\text{\footnotesize{$\rangle$}}$}, then we would have \mbox{$k
> p^d+3$}. But, on the other hand, we have $k \leq p^d+1$ by
Corollary~\ref{Cameron_t=3}, a contradiction. Therefore, we can
argue as in Case (2) to obtain the desired contradiction.

Now, let $p^a=2$. To exclude trivial Steiner \mbox{$3$-designs},
let $v=2^{2d} >4$. For $d=2$ (here $Sp(4,2) \cong S_6$ as
well-known), Corollary~\ref{Cameron_t=3} yields $k \leq 5$. As
$k-2 \nmid v-2$ for $k=5$, it is sufficient by
Lemma~\ref{Comb_t=3}~(d) to consider the case when $k=4$. For $d
>2$, we will show that we can also restrict ourselves to Steiner
quadruple systems. Hence, the claim follows from~\cite{Hu2001}
again. It is easily seen that there are $2^{2d-1}(2^{2d}
- 1)$ hyperbolic pairs in the non-degenerate symplectic space
$V=V(2d,2)$, and by Witt's theorem, $Sp(2d,2)$ is transitive on
these hyperbolic pairs. Let $\{x,y\}$ denote a hyperbolic pair,
and $\mathcal{E}=\text{\footnotesize{$\langle$}} x,y
\text{\footnotesize{$\rangle$}}$ the hyperbolic plane spanned by
$\{x,y\}$. As $\mathcal{E}$ is non-degenerate, we have
the orthogonal decomposition
\[V=\mathcal{E}\perp\mathcal{E}^\perp.\]
Clearly, $Sp(2d,2)_{\{x,y\}}$ stabilizes $\mathcal{E}^\perp$
as a subspace, which implies that \linebreak
\mbox{$Sp(2d,2)_{\{x,y\}} \cong Sp(2d-2,2)$}. As $\left| \Out(Sp(2d,2))
\right|=1$, we have therefore
\[Sp(2d-2,2) \cong Sp(2d,2)_{\{x,y\}} \unlhd
Sp(2d,2)_\mathcal{E}=G_{0,\mathcal{E}}.\] Since $Sp(2d-2,2)$ acts
transitively on the non-zero vectors of the
\mbox{$(2d-2)$}-dimensional symplectic subspace, it is easy to see
that the smallest orbit on \mbox{$V \setminus \mathcal{E}$} under
$G_{0,\mathcal{E}}$ has length at least $2^{2d-2}-1$. If the
unique block $B \in \B$ which is incident with the
\mbox{$3$-subset} $\{0,x,y\}$ contains some point in \mbox{$V
\setminus \mathcal{E}$}, then we would have $k \geq 2^{2d-2}+2$, a
contradiction to Corollary~\ref{Cameron_t=3}. Thus, $B$ lies
completely in $\mathcal{E}$, and with regard to the
flag-transitivity of $G$, we conclude that each block must be
contained in an affine plane. Therefore, we have $k \leq 4$, and
in particular $k=4$ as trivial Steiner \mbox{$3$-designs} are
excluded.

\bigskip
\emph{Case} (4): $G_0 \unrhd G_2(2^a)'$, $d=6a$.
\medskip

We will also show by contradiction that \mbox{$G \leq \Aut(\D)$}
cannot act flag-transitively on any non-trivial Steiner
\mbox{$3$-design} $\D$. First, let $a=1$. Then we have $v=2^6=64$,
and by Corollary~\ref{Cameron_t=3}, it follows that $k \leq 9$.
But, on the other hand, we have $\left| G_2(2)' \right|=2^5 \cdot
3^3 \cdot 7$ and $\left| \Out(G_2(2)') \right| =2$. Thus, in view
of Lemma~\ref{divprop}, we obtain
\[r=\frac{63 \cdot 62}{(k-1)(k-2)}
 \Bigm|  \left| G_0 \right|
\Bigm|  2^6 \cdot 3^3 \cdot 7.\] But this implies that $k-1$ or
$k-2$ is a multiple of $31$, a contradiction.

Now, let $a>1$. As here $G_2(2^a)$ is simple non-Abelian, it is
sufficient to consider $G_0 \unrhd G_2(2^a)$. The permutation
group $G_2(2^a)$ is of \mbox{rank $4$}, and for $0 \neq x \in V$,
the one-point stabilizer $G_2(2^a)_x$ has exactly three orbits
$\mathcal{O}_i$ $(i=1,2,3)$ on \mbox{$V \setminus
\text{\footnotesize{$\langle$}} x
\text{\footnotesize{$\rangle$}}$} of length
$2^{3a}-2^a,2^{5a}-2^{3a},2^{6a}-2^{5a}$ (see,
e.g.,~\cite{Asch1987} or~\cite[Thm.\,3.1]{CaKa1979}). Thus, $G_0$
has exactly three orbits on \mbox{$V \setminus
\text{\footnotesize{$\langle$}} x
\text{\footnotesize{$\rangle$}}$} of length at least $\left|
\mathcal{O}_i \right|.$ Let $y \in \text{\footnotesize{$\langle$}}
x \text{\footnotesize{$\rangle$}}$ be distinct from $0$ and $x$.
Again, we will show that the unique block $B \in \B$ which is
incident with the \mbox{$3$-subset} $\{0,x,y\}$ lies completely in
$\text{\footnotesize{$\langle$}} x
\text{\footnotesize{$\rangle$}}$. If $B$ contains at least one
point of \mbox{$V \setminus \text{\footnotesize{$\langle$}} x
\text{\footnotesize{$\rangle$}}$} in $\mathcal{O}_2$ or
$\mathcal{O}_3$, then we would obtain as above a contradiction to
Corollary~\ref{Cameron_t=3}. Thus, we only have to consider the
case when $B$ contains points of \mbox{$V \setminus
\text{\footnotesize{$\langle$}} x
\text{\footnotesize{$\rangle$}}$} which all lie in
$\mathcal{O}_1$. By~\cite{Asch1987}, the orbit $\mathcal{O}_1$ is
exactly known, and we have
\[\mathcal{O}_1 = x \Delta \setminus
\text{\footnotesize{$\langle$}} x
\text{\footnotesize{$\rangle$}},\] where $x\Delta = \{y \in V \mid
f(x,y,z)=0 \;\,\mbox{for all} \;\, z \in V\}$ with an alternating
trilinear form $f$ on $V$. Then $B$ consists, apart from elements
of $\text{\footnotesize{$\langle$}} x
\text{\footnotesize{$\rangle$}}$, exactly of $\mathcal{O}_1$.
Since $\left| \mathcal{O}_1 \right| \neq 1$, we can choose
$\text{\footnotesize{$\langle$}} \overline{x}
\text{\footnotesize{$\rangle$}} \in x\Delta$ with
$\text{\footnotesize{$\langle$}} \overline{x}
\text{\footnotesize{$\rangle$}} \neq
\text{\footnotesize{$\langle$}} x
\text{\footnotesize{$\rangle$}}$. Let $\overline{y} \in
\text{\footnotesize{$\langle$}} \overline{x}
\text{\footnotesize{$\rangle$}}$ be distinct from $0$ and
$\overline{x}$. Then, for symmetric reasons, the \mbox{$3$-subset}
$\{0,\overline{x},\overline{y}\}$ is also incident with the unique
block $B$. But, on the other hand, we have $\overline{x} \Delta
\neq x \Delta$ for $\text{\footnotesize{$\langle$}} \overline{x}
\text{\footnotesize{$\rangle$}} \neq
\text{\footnotesize{$\langle$}} x
\text{\footnotesize{$\rangle$}}$, a contradiction. Thus, $B$ is
contained completely in $\text{\footnotesize{$\langle$}} x
\text{\footnotesize{$\rangle$}}$, and we may argue as in the cases
above.

\bigskip
\emph{Case} (5): $G_0 \cong A_6$ or $A_7$, $v=2^4$.
\medskip

As $v=2^4$, we have $k \leq 5$ by Corollary~\ref{Cameron_t=3}. If
$k=4$, then applying~\cite{Hu2001} yields the flag-transitive
design described in part (1) of the Main Theorem with assertion
(iii). For $k=5$, we obtain with Lemma~\ref{Comb_t=3}~(d) a
contradiction.

\bigskip
\emph{Cases} (6)-(8).
\medskip

For the existence of non-trivial Steiner \mbox{$3$-designs}, we
have in these cases only a small number of possibilities for $k$
to check, which can easily be ruled out by hand using
Lemma~\ref{Comb_t=3}~(b) and (d), and Corollary~\ref{Cameron_t=3}.

\bigskip


\section{\mbox{Cases with a Group of Automorphisms of Almost Simple
Type}}\label{almost simple type}

Maintaining the same notation, let $\D=(X,\B,I)$ be a
non-trivial Steiner \mbox{$3$-design} with \mbox{$G \leq
\Aut(\D)$} acting flag-transitively on $\D$. We will examine in
this section successively those cases where $G$ is of almost
simple type.

\bigskip
\emph{Case} (1): $N=A_v$, $v \geq 5$. Here, $G$ is
\mbox{$3$-transitive} and does not act on any non-trivial Steiner
\mbox{$3$-design} by~\cite[Thm.\,3]{Kant1985}.

\bigskip

\emph{Case} (2): $N=PSL(d,\tilde{q})$, $d \geq 2$,
$v=\frac{\tilde{q}^d-1}{\tilde{q}-1}$, where $(d,\tilde{q}) \not=
(2,2),(2,3)$.
\medskip

We distinguish two subcases:

\medskip
\indent (i) $N=PSL(2,\tilde{q})$, $v=\tilde{q}+1$.
\smallskip

Let $\tilde{q}=q^e$, $e \geq 1$. Without restriction, we have here
$q^e \geq 5$ as \linebreak $PSL(2,4) \cong PSL(2,5)$, and
$\Aut(N)= P \mathit{\Gamma} L (2,q^e)$. First, we suppose that $G$
is \mbox{$3$-transitive}. In view of~\cite[Thm.\,3]{Kant1985}, we
have then only the \mbox{$3$-$(q^e+1,q+1,1)$} design described in
part (2) of the Main Theorem (without the subcase in brackets)
with \mbox{$PSL(2,q^e) \leq G \leq P \mathit{\Gamma} L(2,q^e)$},
$q \geq 3$, $e \geq 2$. Conversely, flag-transitivity holds as the
$3$-transitivity of $G$ implies that $G_x$ acts block-transitively
on the derived Steiner \mbox{$2$-design} $\D_x$ for any $x \in X$.
Since $PGL(2,q^e)$ is a transitive extension of $AGL(1,q^e)$, it
is easily seen that the derived design at any given point of
\mbox{$GF(q^e) \cup \{\infty\}$} is isomorphic to the
\mbox{$2$-$(q^e,q,1)$} design consisting of the points and lines
of $AG(e,q)$.

Now, we suppose that $G$ is \mbox{$3$-homogeneous} but not
\mbox{$3$-transitive}. Since here $PSL(2,q^e)$ is a transitive
extension of $AG^2L(1,q^e)$ (which is the group of all
permutations of $GF(q^e)$ of the form $x \mapsto a^2x + c$ with
$a,c \in GF(q^e),\,a \neq 0)$, we can deduce
from~\cite{Deletal1986} that the derived design at any given point
is either $AG(e,q)$ with the lines as blocks or the Netto triple
system $N(q^e)$. Thus, part (2) of the Main Theorem holds with the
subcase in brackets or part (3) with \mbox{$PSL(2,q^e) \leq G \leq
P\mathit{\Sigma} L(2,q^e)$} (where, for an odd prime $p$, we
define $P \mathit{\Sigma} L (2,p^a)= PSL(2,p^a) \rtimes
\text{\footnotesize{$\langle$}} \tau_{\alpha}
\text{\footnotesize{$\rangle$}}$ with $\tau_{\alpha} \in$
Sym$(GF(p^a) \cup \{\infty\}) \cong S_v$ of order $a$ induced by
the Frobenius automorphism $\alpha : GF(p^a) \longrightarrow
GF(p^a)$, \mbox{$x \mapsto x^p$}). Conversely, in view of its
\mbox{$3$-homogeneity}, $G$ is also block-transitive. By the
orbit-stabilizer property, we obtain $\left| PSL(2,q^e)_B \right|
= \left| PSL(2,q) \right|$ and in view
of~\cite[Ch.\,12,\,p.\,286]{Dick1901} actually
\[PSL(2,q^e)_B \cong PSL(2,q)\] for any $B \in \B$. Since
$PSL(2,q)$ acts \mbox{$2$-transitively} on $k=q+1$ points, it
follows that in both cases flag-transitivity holds.

Finally, we assume that $G$ is not \mbox{$3$-homogeneous}. As
$PGL(2,q^e)$ is \mbox{$3$-homogeneous}, the unique orbit under
$PGL(2,q^e)$ on the \mbox{$3$-subsets} of $X$ splits under
$PSL(2,q^e)$ in exactly two orbits of equal length. Thus, $G$ has
here exactly two orbits of equal length on the \mbox{$3$-subsets}
of $X$, and by the definition of Steiner \mbox{$3$-designs}, it
follows that $G$ has exactly two orbits (possibly of different
length) on the blocks. Hence, $G\leq \Aut(\D)$ cannot act
block-transitively, and therefore not flag-transitively, on any
non-trivial Steiner \mbox{$3$-design} $\D$.

\medskip
\indent (ii) $N=PSL(d,\tilde{q})$, $d \geq 3$.
\smallskip

We have here $\Aut(N)=P \mathit{\Gamma} L(d,\tilde{q}) \rtimes
\text{\footnotesize{$\langle$}} \iota_\beta
\text{\footnotesize{$\rangle$}}$, where $\iota_\beta$ denotes the
graph automorphism induced by the inverse-transpose map
$\beta:GL(d,\tilde{q}) \longrightarrow GL(d,\tilde{q})$, $x
\mapsto {^t(x^{-1})}$. We will prove by contradiction that $G\leq
\Aut(\D)$ cannot act on any non-trivial Steiner \mbox{$3$-design}
$\D$.

Let us first assume that $d=3$. By the definition of Steiner
\mbox{$3$-designs}, we may choose in the underlying projective
plane $PG(2,\tilde{q})$ three distinct non-collinear points $x,y,z
\in X$, which are incident with a unique block $B \in \B$. We
consider two subcases:

\smallskip
\indent (a) $B$ contains at least one further point of the
triangle through $x, y, z$.

\smallskip
\indent (b) $B$ does not contain any further point of the
triangle.
\smallskip

ad (a): Let $\g$ denote a line of $PG(2,\tilde{q})$. It is
well-known that the translation group $T(\g)$ operates regularly
on the points of $PG(2,\tilde{q})\setminus \g$ and acts trivially
on $\g$. Thus, $T(\g)$ fixes a block $B \in \B$ if three or more
distinct points of $B$ lie on $\g$. Therefore, the block mentioned
in (a) must contain all points of $PG(2,\tilde{q})\setminus \g$,
thus at least $\tilde{q}^2+3$ many. But, these are obviously more
than half of the points of $PG(2,\tilde{q})$, a contradiction to
$k \leq \bigl\lfloor \frac{v}{4} + 2 \bigr\rfloor$ by
Proposition~\ref{Cam}~(a).

ad (b): The pointwise stabilizer of three distinct points in
$SL(3,\tilde{q})$ consists precisely of the diagonal matrices, and
hence has order $(\tilde{q}-1)^2$ (see,
e.g.,~\cite[Ch.\,II,\,Thm.\,7.2\,(b)]{HupI1967}). To this
corresponds in $PSL(3,\tilde{q})$ a subgroup $U$ of order
\[\textstyle{\frac{1}{n}}(\tilde{q}-1)^2 \quad \mbox{with} \quad
n = (3,\tilde{q}-1).\] As $U$ acts semi-regularly outside the
triangle, we obtain $n$ point-orbits of equal length
$\frac{1}{n}(\tilde{q}-1)^2$, since if $U$ fixes some further
point outside the triangle, then $U$ would fix some non-degenerate
quadrangle, and so would be the identity, a contradiction. Thus,
we get
\[k \geq 3+\textstyle{\frac{1}{n}}(\tilde{q}-1)^2. \]
On the other hand, we know that the block mentioned in (b) is an
arc, and therefore contains at most $\tilde{q}+1$ points for
$\tilde{q}$ odd or $\tilde{q}+2$ points for $\tilde{q}$
even (see, e.g.,~\cite[Ch.\,3.2,\,Thm.\,24]{Demb1968}). Only for
$\tilde{q}=2$ and $4$ both conditions are fulfilled. But, with
regard to Lemma~\ref{Comb_t=3}~(d), there exist no non-trivial
\mbox{$3$-$(7,k,1)$ designs} and \mbox{$3$-$(21,k,1)$ designs}.
Therefore, for $d=3$ we have shown that $G$ cannot act on any
non-trivial \mbox{$3$-$(\tilde{q}^2+\tilde{q}+1,k,1)$}
design.

Now, we consider the case when $d > 3$. Via induction over $d$, we
will verify that $G\leq \Aut(\D)$ cannot act on any non-trivial
Steiner \mbox{$3$-design} $\D$. For this, let us assume that there
is a counter-example with $d$ minimal. Without restriction, we can
choose three distinct points $x,y,z$ from a hyperplane $\h$ of
\mbox{$PG(d-1,\tilde{q})$}. First, we show that the unique block
$B \in \B$ which is incident with the \mbox{$3$-subset}
$\{x,y,z\}$ is contained completely in $\h$. Analogously as above,
the translation group $T(\h)$ acts regularly on the points of
\mbox{$PG(d-1,\tilde{q}) \setminus \h$}, but trivially on $\h$. If
$B$ contains at least one point outside $\h$, then it would
already contain all points of $PG(d-1,\tilde{q}) \setminus \h$,
thus at least $\tilde{q}^{d-1}+3$ many. However, as
\[v=\frac{\tilde{q}^d-1}{\tilde{q}-1} < 2 \tilde{q}^{d-1}
\Longleftrightarrow \tilde{q}^d-1 <
2(\tilde{q}^d-\tilde{q}^{d-1})\Longleftrightarrow 2
\tilde{q}^{d-1}-1 < \tilde{q}^d,\] these are more than half of
the points of $PG(d-1,\tilde{q})$, the same contradiction as
above. Thus, $\h$ induces a
\[3\mbox{-}(\textstyle{\frac{\tilde{q}^{d-1}-1}{\tilde{q}-1}},k,1) \;
\mbox{design},\] on which $G$ containing $PSL(d-1,\tilde{q})$
as simple normal subgroup operates. Inductively, we obtain the
minimal counter-example for $d=3$. But, as we have shown above, $G$
with $PSL(3,\tilde{q})$ as simple normal subgroup cannot act on
any non-trivial \mbox{$3$-$(\tilde{q}^2+\tilde{q}+1,k,1)$}
design, and the assertion follows.

\bigskip
\emph{Case} (3): $N=PSU(3,q^2)$, $v=q^3+1$, $q=p^e>2$.
\medskip

Here $\Aut(N)= P \mathit{\Gamma} U(3,q^2)$, and $\left| G \right|
= (q^3+1)q^3 \frac{(q^2-1)}{n}a$ with $n=(3,q+1)$ and $a \mid 2
ne$. Thus, from Remark~\ref{equa_t=3}, we obtain
\begin{equation}\label{eq-3}
q^2+q+1 = (k-1)(k-2) \frac{q+1}{n}\frac{a}{\left|G_{xB}\right|}
\quad \mbox{if} \;\, x \in B.
\end{equation}
We will show by contradiction that $G\leq \Aut(\D)$ cannot act
flag-transitively on any non-trivial Steiner \mbox{$3$-design} $\D$.

Let $\{v_1,v_2,v_3\}$ be a basis of the non-degenerate hermitian
vector space $V=V(3,q^2)$ with
\[(v_2,v_2)=(v_1,v_3)=1,\;(v_1,v_1)=(v_3,v_3)=(v_1,v_2)=(v_2,v_3)=0.\]
For $v =\sum_{i=1}^{3} a_i v_i$ and $w =\sum_{i=1}^{3} b_i v_i$
($a_i,\,b_i \in GF(q^2)$), we have then
\[(v,w)=a_1b_3^\tau + a_2 b_2^\tau + a_3b_1^\tau,\]
where $\tau$ denotes the unique involutory automorphism $GF(q^2)
\longrightarrow GF(q^2)$, \mbox{$x \mapsto x^q$}. We deduce from
the proof of~\cite[Ch.\,II,\,Thm.\,10.12]{HupI1967} that the
cyclic group
\[\left\{ \begin{pmatrix}
  c &        & \\
    & c^{-2} & \\
    &        & c
\end{pmatrix} \Biggm| c^{-2} \neq c,\,c \in GF(q^2)^* \right\} \]
of linear transformations on $V$ induces a group $U$ of
dilatations of order $\frac{q+1}{n}$ on the associated projective
space $PG(2,q^2)$ with axis the non-absolute line $\g$ consisting
of the absolute points
$\text{{$\langle$}}(1,0,0)\text{{$\rangle$}}$,
$\text{{$\langle$}}(0,0,1)\text{{$\rangle$}}$ and
$\text{{$\langle$}}(a_1,0,a_3)\text{{$\rangle$}}$ with
\[a_1a_3^{\tau} + a_1^{\tau} a_3 = \mbox{Tr}(a_1 a_3^{\tau})=0\]
(where $\mbox{Tr}$ denotes the trace map $GF(q^2) \longrightarrow
GF(q)$, $x \mapsto x + x^q$)
and as center the pole of the axis, i.e. the non-absolute point
$\text{{$\langle$}}(0,1,0)\text{{$\rangle$}}$.

As it is customary (see, e.g.,~\cite[p.\,87]{Baer1946}), we call
in the following non-absolute lines $\g$ and $\h$ \emph{perpendicular}
if $\g$ passes through the pole of $\h$ and $\h$ passes, therefore,
through the pole of $\g$.

By the definition of Steiner \mbox{$3$-designs}, we may choose
three distinct absolute points on $\g$, which are incident with a
unique block $B \in \B$. Let us first assume that $B$ contains
absolute points outside $\g$ which are all on $\h$. It is clear
that $U$ fixes each point of $\g$, and hence in particular $B$.
Furthermore, $\h$ intersects $\g$ in a non-absolute point $x$
(see, e.g.,~\cite[p.\,88]{Baer1946}). As $U$ acts outside $x$
semi-regularly on $\h$, we conclude that all point-orbits have
length $\frac{q+1}{n}$. If we choose now three distinct absolute
points on $\h$, then they are also incident with the unique block
$B$. Thus, by the same arguments, $U$ fixes each point of $\h$ and
acts outside $x$ semi-regularly on $\g$. Therefore, we have
\[k= (n_1 + n_2) \frac{q+1}{n}\]
with $n_1,n_2 \in \{1,2,3\}$. If $n=1$, then obviously $k=2(q+1)$,
which is impossible in view of Lemma~\ref{Comb_t=3}~(d). Thus, $n
\neq 1$. For $n_1+n_2=3$, it follows from equation~(\ref{eq-3})
that $q^2+q+1 \bigm| (q-1) \frac{a}{n} < q^2-q$, which is clearly
not possible. In each of the other cases, polynomial division with
remainder gives a contradiction to Lemma~\ref{Comb_t=3}~(d).

Now, we assume that $B$ contains absolute points outside $\g$
which are not all on $\h$. By applying the same arguments as
above, we obtain additionally a lattice of points such that
\[k= n_1n_2\left( \frac{q+1}{n} \right)^2 +
(n_1 + n_2)\frac{q+1}{n}\] with $n_1, n_2$ as above,
which clearly contradicts Corollary~\ref{Cameron_t=3}.

Hence, we have shown that $B$ is completely contained in $\g$.
Thus, in view of the flag-transitivity of $G$, each block is
contained in a non-absolute line. But, by the definition of
Steiner \mbox{$3$-designs}, any three non-collinear absolute
points must also be incident with a unique block, a contradiction.

\bigskip
\emph{Case} (4): $N=Sz(q)$, $v=q^2+1$, $q=2^{2e+1}>2$.
\medskip

We have $\Aut(N)= Sz(q) \rtimes \text{\footnotesize{$\langle$}}
\alpha \text{\footnotesize{$\rangle$}}$, where $\alpha$ denotes
the Frobenius automorphism $GF(q) \longrightarrow GF(q)$, \mbox{$x
\mapsto x^2$}. Thus, by Dedekind's law, $G = Sz(q) \rtimes (G \cap
\text{\footnotesize{$\langle$}} \alpha
\text{\footnotesize{$\rangle$}})$, and $\left| G \right| =
(q^2+1)q^2(q-1)a$ with $a \mid 2e+1$. It follows from
Remark~\ref{equa_t=3} that
\[ q+1 = (k-1)(k-2) \frac{a}{\left|G_{xB}\right|} \quad
\mbox{if} \;\, x \in B.\] We will prove by contradiction that
\mbox{$G \leq \Aut(\D)$} cannot act flag-transitively on any
non-trivial Steiner \mbox{$3$-design} $\D$.

Let us first remark that we only have one class of involutions in
$G$. Hence, every involution has exactly one fixed point, which
lies in an appropriate block. Therefore, by the flag-transitivity
of $G$, there exists for every $B \in \B$ always an involution
$\tau \in G_{xB} \cap Sz(q)$ with $x \in B$, and $B$ can be
regarded as the orbit of fixed points of involutions in $G_B \cap
Sz(q)$.

Since $G$ is block-transitive, we can restrict ourselves to
consider the unique block $B \in \B$ which is incident with the
\mbox{$3$-subset} $\{0,1,\infty\}$ of $X$. As every non-identity
element of $Sz(q)$ fixes at most two distinct points, we have
$\Aut(N)_{0,1,\infty} = \text{\footnotesize{$\langle$}} \alpha
\text{\footnotesize{$\rangle$}}$, and thus $G \cap
\text{\footnotesize{$\langle$}} \alpha
\text{\footnotesize{$\rangle$}} \leq G_{0B}$ by the definition of
Steiner \mbox{$3$-designs}. Setting $u=\frac{\left| G_{0B}
\right|}{a}$, we next show that $u=2$ or $4$. For the list of
subgroups of $Sz(q)$, we refer to~\cite[Thm.\,9]{Suz1962}. First,
let $G_B \cap Sz(q)$ be isomorphic to $Sz(\overline{q})$ for some
$\overline{q} \geq 8$ such that $\overline{q}^m = q$, $m \geq 1$.
As $B$ can be regarded as the orbit of fixed points of involutions
in $G_B \cap Sz(q)$, it follows that $k = \overline{q}^2+1$.
Clearly, $m > 1$ (otherwise, $k=q^2+1$, a contradiction to
Corollary~\ref{Cameron_t=3}). Thus, we have
\[q+1=\overline{q}^2(\overline{q}^2-1) \frac{a}{\left| G_{0B} \right|}\,.\]
As $q>8$, Zsigmondy's theorem yields the existence of a
$2$-primitive prime divisor $\overline{r}$ with $\overline{r}
\perp 2^{2(2e+1)}-1$. Then
\[\overline{r}  \bigm|  q+1 = \overline{q}^2(\overline{q}^2-1)
\frac{a}{\left|G_{0B}\right|}\,.\] But
now~\cite[Thm.\,3.5\,(ii)]{Her1974} yields
$(\overline{r},\overline{q})=1$ and $\overline{r} > a$ since
$\overline{r} \equiv 1$ (mod $(2e+1)$). Therefore, we conclude that
$\overline{q}=q$, a contradiction.\\
Let $G_B \cap Sz(q)$ be conjugated to a subgroup of $Sz(q)_x$ $(x
\in X)$. By the transitivity of $G$, we can choose $x$ as fixed
point of an involution. Thus, $x \in B$ by the remark above,
contrary to the fact that $x \notin B$ by the flag-transitivity of $G$.\\
\noindent Let $G_B \cap Sz(q)$ be conjugated to a subgroup of $U$
with \mbox{$\left| U \right|=4(q \pm l +1)$,} where $l^2=2q$. Then
$\left| O_{p'} (U) \right|= q \pm l +1$, and $O_{p'} (U)$ operates
fixed-point-\linebreak freely on $X$ since $(q \pm l +1,q)=1$ and $(q \pm l
+1,q^2- 1)=1$. Thus \mbox{$(G_{0B} \cap Sz(q)) \cap O_{p'} (U)
=1$}, and therefore $\left|G_{0B} \cap Sz(q)
\right| \leq 4$.\\
Let $G_B \cap Sz(q)$ be conjugated to a subgroup of $U$ with
$\left| U \right| = 2(q-1)$. Then $\left| O_{p'} (U) \right|=
q-1$, and $O_{p'} (U)$ has two distinct fixed points in $X$. As
$O_{p'} (U)$ contains no involutions, these fixed points cannot
lie in $B$ by the remark above. Hence \mbox{$(G_{0B} \cap Sz(q))
\cap O_{p'} (U) =1$}, and thus $\left|
G_{0B} \cap Sz(q) \right| \leq 2$.\\
Since $\left|G_{0B} \cap Sz(q) \right| \equiv 0$ (mod $2)$, we
have therefore \[\left|G_{0B} \cap Sz(q) \right| = 2 \; \mbox{or}
\; 4.\] As $G \cap \text{\footnotesize{$\langle$}} \alpha
\text{\footnotesize{$\rangle$}} \leq G_{0B}$, and clearly $(G_B
\cap Sz(q)) \cap (G \cap \text{\footnotesize{$\langle$}} \alpha
\text{\footnotesize{$\rangle$}}) =1$, we conclude that
\[u= 2 \; \mbox{or} \; 4.\]

Finally, our equation
\[u(q+1)=(k-1)(k-2)\]
yields for $u=2$ that
\[2^{2e+2}= k (k-3),\]
which is clearly impossible since $e \geq 1$. For $u=4$, we obtain
\begin{equation}\label{eq-4}
2^{2e+3}=k^2-3k-2.
\end{equation}
By setting $x=2k-3$ and $n=2e+5$ this becomes the well-known
generalized Ramanujan-Nagell equation
\[x^2-17 = 2^n,\]
which has exactly the four solutions
$(x,n)=(5,3),\,(7,5),\,(9,6),\,(23,9)$ (see,
e.g.,~\cite[Thm.\,3]{Beuk1981}). As we have $e \geq 1$, it follows
that $(e,k)=(2,13)$ is the only solution of equation~(\ref{eq-4}).
But, by Lemma~\ref{Comb_t=3}~(b), this is impossible, which
verifies the claim.

\bigskip
\emph{Case} (5): $N=Re(q)$, $v=q^3+1$, $q=3^{2e+1}>3$.
\medskip

Here $\Aut(N)= Re(q) \rtimes \text{\footnotesize{$\langle$}}
\alpha \text{\footnotesize{$\rangle$}}$, where $\alpha$ denotes
the Frobenius automorphism $GF(q) \longrightarrow GF(q)$, \mbox{$x
\mapsto x^3$}. Thus, by Dedekind's law, $G = Re(q) \rtimes (G \cap
\text{\footnotesize{$\langle$}} \alpha
\text{\footnotesize{$\rangle$}})$, and $\left| G \right| =
(q^3+1)q^3(q-1)a$ with $a \mid 2e+1$. From Remark~\ref{equa_t=3},
we hence obtain
\begin{equation}\label{eq-5}
q^2+q+1 = (k-1)(k-2) \frac{a}{\left|G_{xB}\right|} \quad
\mbox{if} \;\, x \in B.
\end{equation}
We will also prove by contradiction that $G \leq \Aut(\D)$
cannot act flag-transitively on any non-trivial Steiner
\mbox{$3$-design} $\D$.

We remark that we only have one class of involutions in $G$. Thus,
every involution fixes at least three distinct points, each of
which lies in an appropriate block. Therefore, by the
flag-transitivity of $G$, there exists for every $B \in \B$ always
an involution $\tau \in G_{xB} \cap Re(q)$ with $x \in B$.

We show furthermore that $9 \bigm| \left|G_B \cap Re(q) \right|$.
Let $P$ be a Sylow $3$-subgroup of $Re(q).$ According
to~\cite{Ward1966}, $P$ contains a normal elementary Abelian
subgroup $\overline{P}$ of order $q^2$ containing $Z(P)$. Thus,
there exist subgroups $U_1$, $U_2$ of $\overline{P}$ of order $3$
with $U_1 \leq Z(P)$, $U_2 \nleq Z(P)$. As the stabilizer of three
distinct points in $Re(q)$ has order $2$, we have
$\mbox{Fix}_X(U_1) = \mbox{Fix}_X(U_2) = \{x\}$ for some $x \in
X$. Hence, if $U_1$ and $U_2$ are conjugated in $Re(q)$, then they
are already conjugated in $Re(q)_x$. But, as $Z(P)$ is a
characteristic subgroup of $Re(q)_x$, this is impossible.
Therefore, we have at least two distinct classes of subgroups of
order $3$ in $Re(q)$, and the assertion follows by the definition
of Steiner \mbox{$3$-designs}.

Because of the block-transitivity of $G$, we can restrict
ourselves to consider the unique block $B \in \B$ which is
incident with the \mbox{$3$-subset} $\{0,1,\infty\}$ of $X$.
Clearly, $\text{\footnotesize{$\langle$}} \alpha
\text{\footnotesize{$\rangle$}} \leq \Aut(N)_{0,1,\infty}$, and
hence \mbox{$G \cap \text{\footnotesize{$\langle$}} \alpha
\text{\footnotesize{$\rangle$}} \leq G_{0B}$} by the definition of
Steiner \mbox{$3$-designs}. Furthermore, obviously $(G_B \cap
Re(q)) \cap (G \cap \text{\footnotesize{$\langle$}} \alpha
\text{\footnotesize{$\rangle$}}) =1$. Therefore, as $G_B$ acts
transitively on the points of $B$, Dedekind's law yields
\begin{equation}\label{eq-6}
k=\left|0^{G_B}\right| = \big[G_B : G_{0B} \big] =
\big[ G_B \cap Re(q) : G_{0B} \cap Re(q) \big].
\end{equation}
Thus, $G_B \cap Re(q)$ acts also transitively on the points of
$B$.

In the following, we will examine the list of subgroups of $Re(q)$
(cf.~\cite{Ward1966}). As $9$ divides the order of \mbox{$G_B \cap
Re(q)$}, clearly \mbox{$G_B \cap Re(q)$} cannot be conjugated to a
subgroup of the normalizer of a Sylow \mbox{$2$-subgroup} of
$Re(q)$ of order $8 \cdot 7 \cdot 3$. By the same argument, $G_B
\cap Re(q)$ cannot be conjugated to a
subgroup of $U$ with $\left| U \right| = 6(q+1 \pm 3l)$, where $l=3^e$.\\
Let $G_B \cap Re(q)$ be isomorphic to $Re(\overline{q})$ for some
$\overline{q} \geq 27$ such that $\overline{q}^m = q$, $m \geq 1$.
Let $\overline{X} \subseteq X$ with $\left| \overline{X}
\right|=\overline{q}^3+1$. We first show that only involutions may
have fixed points in $X \setminus \overline{X}$. Let $g \in G$
with $o(g)=s$, where $s \neq 2$ is a prime. If $s \mid
\overline{q}-1$, then $g$ has two distinct fixed points in
$\overline{X}$, and none in $X \setminus \overline{X}$ since the
stabilizer of three distinct points in $Re(q)$ has order $2$. For
$s=3$, clearly $g$ has exactly one fixed point, which lies in
$\overline{X}$. If $s \mid \overline{q}+1$, we show that $g$ has
no fixed point in $X$. Obviously, $g$ has no fixed point in
$\overline{X}$. As $3 \nmid \overline{q}+1$, we assume that $g$
has two distinct fixed points in $X \setminus \overline{X}$. But,
as
\[q^3-\overline{q}^3 = \bigg( \sum_{i=0}^{3m-1}(-1)^i \frac{q^3}
{\overline{q}^{i+1}} - \overline{q}^2 + \overline{q}-1 \bigg)
\bigg(\overline{q}+1 \bigg),\] and hence
$(q^3-\overline{q}^3-2,\overline{q}+1)=(2,\overline{q}+1)=2$, this
is impossible. If \mbox{$s \mid \overline{q}+1 \pm 3\overline{l}$},
we show again that $g$ has no fixed point in $X$. As
$\overline{q}^3+1=(\overline{q}+1+3\overline{l})(\overline{q}+1-3
\overline{l})(\overline{q}+1)$, it is obvious that $g$ has no
fixed point in $\overline{X}$. Since $3 \nmid \overline{q}+1 \pm 3
\overline{l}$, we assume in both cases that $g$ has two distinct
fixed points in $X \setminus \overline{X}$. But, as
$(\overline{q}+1+3\overline{l})(\overline{q}+1-3\overline{l})=
\overline{q}^2-\overline{q}+1$,
and
\[q^3- \overline{q}^3 = \bigg( \sum_{i=0}^{m-1} \sum_{j=0}^{1}(-1)^{2+3i}
\frac{q^3}{\overline{q}^{2+3i+j}} - \overline{q}-1 \bigg)
\bigg(\overline{q}^2 -\overline{q}+1 \bigg),\] we have $(q^3 -
\overline{q}^3-2,\overline{q}^2 -
\overline{q}+1)=(2,\overline{q}^2 - \overline{q}+1)=2$, a
contradiction.\\
As $G_B \cap Re(q)$ acts transitively on the points of $B$, we
have $B \subseteq \overline{X}$ or $B \subseteq X \setminus
\overline{X}$. In the first case, equation~(\ref{eq-6}) yields
\[k= \overline{q}^3+1,\] while in the second
\[k=\frac{(\overline{q}^3+1)\overline{q}^3(\overline{q}-1)}{n}\,,\]
where $n$ is a power of $2$, and $n \leq 8$ as the order of
$Re(q)$ is divisible by $8$ but not by $16$.\\
We will prove now that none of these values of $k$ is possible. We
assume first that $k=\overline{q}^3+1$. Clearly, $m
> 1$ (otherwise, $k=q^3+1$, a contradiction to
Corollary~\ref{Cameron_t=3}). Thus, we have
\[q^2+q+1=\overline{q}^3(\overline{q}^3-1) \frac{a}{\left| G_{0B} \right|}.\]
Zsigmondy's theorem yields the existence of a $3$-primitive prime
divisor $\overline{r}$ with $\overline{r} \perp 3^{3(2e+1)}-1$.
Then
\[\overline{r}  \bigm|  q^2+q+1 = \overline{q}^3(\overline{q}^3-1)
\frac{a}{\left|G_{0B}\right|}.\] But
now~\cite[Thm.\,3.5\,(ii)]{Her1974} yields
$(\overline{r},\overline{q})=1$ and $\overline{r} > a$ since
$\overline{r} \equiv 1$ (mod $(2e+1)$). Therefore, we have
$\overline{q}=q$, a contradiction.\\
Now, we assume that
$k=\frac{(\overline{q}^3+1)\overline{q}^3(\overline{q}-1)}{n}$.
Then
\[\left| G_{0B} \right| = \frac{\left| G_B \cap Re(q)\right|
\overline{a}}{k}=n \overline{a},\] where $\overline{a} \mid a$.
Here, $n<4$ since otherwise $(k-1)(k-2) \equiv 0$ (mod $4$) by
equation~(\ref{eq-5}) and, by applying Lemma~\ref{Comb_t=3}~(c),
this would imply that $q^3-1$ is divisible by $4$, which is
impossible since $q-1 \equiv 2$ (mod $8$) in $Re(q)$.\\
Thus, we may assume that $n=2$. Polynomial division with remainder
gives
\begin{eqnarray*}
q^3-1 &=&
\bigg(\sum_{i=0}^{\overline{m}}\frac{2^{2i+1}q^3}{\big((\overline{q}^3+1)\overline{q}^3
(\overline{q}-1)\big)^{i+1}}\bigg)
\bigg(\frac{(\overline{q}^3+1)\overline{q}^3(\overline{q}-1)}{2} - 2 \bigg)\\
& & +
\frac{2^{2\overline{m}+2}q^3}{\big((\overline{q}^3+1)\overline{q}^3
(\overline{q}-1)\big)^{\overline{m}+1}}-1
\end{eqnarray*} for a suitable $\overline{m} \in \N$ (such that
\[\mbox{deg}\bigg(\frac{2^{2\overline{m}+2}q^3}
{\big((\overline{q}^3+1)\overline{q}^3(\overline{q}-1)\big)^{\overline{m}+1}}-1
\bigg) < \mbox{deg} \bigg(\frac{(\overline{q}^3+1)
\overline{q}^3(\overline{q}-1)}{2} - 2 \bigg)\] as is well-known).
As $8 \bigm|\left| Re(\overline{q}) \right|$, clearly
\mbox{$\big((\overline{q}^3+1)\overline{q}^3(\overline{q}-1)\big)^{\overline{m}+1}$}
is divisible by $2^{3(\overline{m}+1)}$. Thus
$\frac{2^{2\overline{m}+2}q^3}{\big((\overline{q}^3+1)\overline{q}^3
(\overline{q}-1)\big)^{\overline{m}+1}}\neq 1$,
yielding a contradiction to Lemma~\ref{Comb_t=3}~(d).\\
Let $G_B \cap Re(q)$ be conjugated to a subgroup of $Re(q)_x$ $(x
\in X)$. By the transitivity of $G$, we can choose $x$ as fixed
point of an involution. Thus, $x \in B$ for an appropriate block
$B \in \B$ by the remark above, contrary to the fact that $x
\notin B$ by the flag-transitivity of $G$.\\
Let $G_B \cap Re(q)$ be conjugated to a subgroup of $PSL(2,q)
\times \text{\footnotesize{$\langle$}}\tau
\text{\footnotesize{$\rangle$}}$, where $\tau$ denotes any
involution in $Re(q)$. By the remark above, we can choose $\tau$
such that $0$ is a fixed point under $\tau$. As $9$ must be a
divisor of the order of \mbox{$G_B \cap Re(q)$}, we can restrict
ourselves to the examination of the following cases
(cf.~\cite[Ch.\,12,\,p.\,285f.]{Dick1901}
or~\cite[Ch.\,II,\,Thm.\,8.27]{HupI1967}):

\begin{enumerate}
\item[(i)] $G_B \cap Re(q)$ is conjugated to $PSL(2,\overline{q})$
or $PSL(2,\overline{q})\times \text{\footnotesize{$\langle$}}\tau
\text{\footnotesize{$\rangle$}}$ for some $\overline{q} \geq 27$
such that $\overline{q}^m = q$, $m \geq 1$.\\
Let $\overline{X} \subseteq X$ with $\left| \overline{X}
\right|=\overline{q}+1$. First, we show again that only
involutions may have fixed points in $X \setminus \overline{X}$.
Let $g \in G$ with $o(g)=s$, where $s \neq 2$ is a prime. If $s
\mid \overline{q}-1$, then $g$ has two distinct fixed points in
$\overline{X}$ and none in $X \setminus \overline{X}$. For $s=3$,
clearly $g$ has exactly one fixed point, which lies in
$\overline{X}$. If $s \mid \overline{q}+1$, we show that $g$ has
no fixed point in $X$. Obviously, $g$ has no fixed point in
$\overline{X}$. As $3 \nmid \overline{q}+1$, we assume that $g$
has two distinct fixed points in $X \setminus \overline{X}$. But,
as
\[q^3-\overline{q} = \bigg( \sum_{i=0}^{3\overline{m}-1}(-1)^i
\frac{q^3} {\overline{q}^{i+1}}-1\bigg) \bigg(\overline{q}+1
\bigg),\] and hence
$(q^3-\overline{q}-2,\overline{q}+1)=(2,\overline{q}+1)=2$, this
is impossible.\\
Again, we have $B \subseteq \overline{X}$ or $B \subseteq X
\setminus \overline{X}$. With equation~(\ref{eq-6}), we obtain
\[k= \overline{q}+1,\]
in the first case, while in the second
\[k=\frac{\overline{q}(\overline{q}^2-1)}{n}\,,\] where $n$ is a
power of $2$, and $n \leq 8$ again.\\
We will prove now that none of the values of $k$ is possible. We
assume first that $k=\overline{q}+1$. Then
\[q^2+q+1 \bigm| \overline{q}(\overline{q}-1)a\]
by equation~(\ref{eq-5}). Since $(q^2+q+1,\overline{q})=1$ and
$(q^2+q+1,\overline{q}-1)=(3,\overline{q}-1)=1$, this is
equivalent to
\[q^2+q+1 \bigm| a,\] which is impossible as clearly $a \leq q$.
Now, we assume that $k=\frac{\overline{q}(\overline{q}^2-1)}{n}$.
Then
\[\left| G_{0B} \right| = \frac{\left| G_B \cap Re(q)\right|
\overline{a}}{k}=\frac{n\overline{a}}{2} \quad \mbox{or} \quad n
\overline{a},\] where $\overline{a} \mid a$. Considering the first
yields
\[(q^2+q+1)\frac{n}{2}=(k-1)(k-2)\frac{a}{\overline{a}}\] by
equation~(\ref{eq-5}). Clearly, $n=2$ is impossible. If $n=4$,
then $k=\frac{\overline{q}(\overline{q}^2-1)}{4}$ is divisible by
$2$ but not by $4$. Thus, $4$ is a divisor of $k-2$, but not of
the left side. For $n=8$, we have $(k-1)(k-2) \equiv
0$ (mod $4$), which is not possible as we have seen above.\\
Now, we assume that $\left| G_{0B} \right| = n \overline{a}$.
Here, $n < 4$ again. For $n=2$, we have
$k=\frac{\overline{q}^3-\overline{q}}{2}$. Then, polynomial
division with remainder gives
\[q^3-1=\bigg(\sum_{i=0}^{\overline{m}}\frac{2^{2i+1}q^3}
{(\overline{q}^3-\overline{q})^{i+1}}\bigg)
\bigg(\frac{\overline{q}^3-\overline{q}}{2} - 2 \bigg) +
\frac{2^{2\overline{m}+2}q^3}{(\overline{q}^3-\overline{q})^
{\overline{m}+1}}-1\] for a suitable $\overline{m} \in \N$. As
$(\overline{q}^2-1)^{\overline{m}+1}$ is divisible by
$2^{3(\overline{m}+1)}$, clearly
$\frac{2^{2\overline{m}+2}q^3}{(\overline{q}^3-\overline{q})^
{\overline{m}+1}}\neq 1$, yielding a contradiction to
Lemma~\ref{Comb_t=3}~(d).

\item[(ii)] $G_B \cap Re(q)$ is conjugated to $U$ or $U \times
\text{\footnotesize{$\langle$}}\tau
\text{\footnotesize{$\rangle$}}$, where $U$ is an elementary
Abelian subgroup of order $\overline{q} \mid q$ of $PSL(2,q)$.\\
Let $\overline{X}\subseteq X$ with $\left| \overline{X}\right| = q
+1$. Clearly, $U$ operates regularly on $\overline{q}$ points, and
each non-identity element of $U$ has $\infty$ as only fixed point
in $\overline{X}$ and none in $X \setminus \overline{X}$.\\ As $2
\nmid \left| U \right|$, it follows that $k= \overline{q}$ in both
of the cases $B \subseteq \overline{X}$ and $B \subseteq X
\setminus \overline{X}$. But, polynomial division with remainder
gives
\[q^3-1= \bigg( \sum_{i=0}^{\overline{m}} \frac{2^i q^3}
{\overline{q}^{i+1}}\bigg) \bigg(\overline{q}-2 \bigg) +
\frac{2^{\overline{m}+1} q^3}{\overline{q}^{\overline{m}+1}}-1\]
for a suitable $\overline{m} \in \N$. As clearly
$\frac{2^{\overline{m}+1} q^3}{\overline{q}^{\overline{m}+1}} \neq
1$, this leads to a contradiction to Lemma~\ref{Comb_t=3}~(d)
again.

\item[(iii)] $G_B \cap Re(q)$ is conjugated to $U$ or $U \times
\text{\footnotesize{$\langle$}}\tau
\text{\footnotesize{$\rangle$}}$, where $U$ is a semi-direct
product of an elementary Abelian subgroup of order $\overline{q}
\mid q$ with a cyclic subgroup of order $c$ of $PSL(2,q)$ with $c
\mid \overline{q}-1$ and $c \mid q-1$.\\
Let $\overline{X} \subseteq X$ with $\left| \overline{X}
\right|=q+1$. Again, we show that only involutions may have fixed
points in $X \setminus \overline{X}$. Let $g \in G$ with $o(g)=s$,
where $s \neq 2$ is a prime. If $s=3$, then $g$ has exactly one
fixed point, which lies in $\overline{X}$. If $s \mid c$, then $g$
has
exactly two distinct fixed points, which lie in $\overline{X}$.\\
For $B \subseteq \overline{X}$, we deduce that $k=\overline{q}$ or
$\overline{q} c$, and for $B \subseteq X \setminus \overline{X}$
that $k= \frac{\overline{q}c}{n}$ with $n \leq 2$ since $q-1
\equiv 2$ (mod $8$). Again, we will prove that none of the values
of $k$ is possible. For $k=\overline{q}$, we have already shown
that this is impossible. We assume next that $k=\overline{q} c$.
If $2 \mid c$, then $k$ is divisible by $2$ but not by $4$.
Therefore, $k-2 \equiv 0$ (mod $4$), and hence \mbox{$q^3-1 \equiv
0$ (mod $4$)} by Lemma~\ref{Comb_t=3}~(d), which is impossible as
we have already seen. For $2 \nmid c$, polynomial division with
remainder gives
\[q^3-1=\bigg(\sum_{i=0}^{\overline{m}}\frac{2^i q^3}
{(\overline{q}c)^{i+1}}\bigg) \bigg(\overline{q}c - 2\bigg) +
\frac{2^{\overline{m}+1}
q^3}{(\overline{q}c)^{\overline{m}+1}}-1\] for a suitable
$\overline{m} \in \N$. But obviously $\frac{2^{\overline{m}+1}
q^3}{(\overline{q}c)^{\overline{m}+1}} \neq 1$, which
leads to the same contradiction as before.\\
Now, we assume that $k=\frac{\overline{q} c}{n}$. Then
\[\left| G_{0B} \right| = \frac{\left| G_B \cap Re(q)\right|
\overline{a}}{k}=n \overline{a} \quad \mbox{or} \quad 2n
\overline{a},\] where $\overline{a} \mid a$. When considering the
first possibility, clearly equation~(\ref{eq-5}) rules out the
case $n=1$. So, we assume that $n=2$. Hence $k= \frac{\overline{q}
c}{2}$, but polynomial division with remainder yields
\[q^3-1=\bigg(\sum_{i=0}^{\overline{m}}\frac{2^{2i+1}q^3}
{(\overline{q}c)^{i+1}}\bigg) \bigg(\frac{\overline{q}c}{2} - 2
\bigg) +
\frac{2^{2\overline{m}+2}q^3}{(\overline{q}c)^{\overline{m}+1}}-1\]
for a suitable $\overline{m} \in \N$. But since $c \mid q-1$, the
largest possible power of $2$ that is contained in
$c^{\overline{m}+1}$ is $2^{\overline{m}+1}$. Thus
$\frac{2^{2\overline{m}+2}q^3}{(\overline{q}c)^{\overline{m}+1}}\neq
1$, the same contradiction as above.\\
Now, we assume that $\left| G_{0B} \right| = 2n \overline{a}$. For
$n=1$, we get $k=\overline{q}c$, which is not possible as shown
above. The case $n=2$ is ruled out by equation~(\ref{eq-5}) since
$(k-1)(k-2)$ is not divisible by $4$ as we already know.
\end{enumerate}
This completes the list of subgroups that we have to examine, and
the claim is established.

\bigskip
\emph{Case} (6): $N=Sp(2d,2)$, $d \geq 3$, $v = 2^{2d-1} \pm
2^{d-1}$.
\medskip

As here \mbox{$\left|\Out(N) \right|=1$}, we have $N=G$. Let $X^+$
(respectively $X^-$) denote the set of points on which $G$
operates. It is well-known that $G_x$ acts on \linebreak $X^{\pm}
\setminus \{x\}$ as $O^\pm(2d,2)$ does in its usual rank $3$
representation on singular points of the underlying orthogonal
space. Thus, $G_{xy}$ has two orbits on $X^{\pm}\setminus \{x,y\}$
of length $2(2^{d-1} \mp 1)(2^{d-2} \pm 1)$ and $2^{2d-2}$ (see,
e.g.,~\cite[p.\,69]{Kant1985}). We will show by contradiction that
\mbox{$G \leq \Aut(\D)$} cannot act flag-transitively on any
non-trivial Steiner \mbox{$3$-design} $\D$.

Let $z \in X^{\pm} \setminus \{x,y\}$. Then, in both cases, the
\mbox{$3$-subset} $\{x,y,z\}$ is incident with a unique block $B
\in \B$. By Remark~\ref{equa_t=3}, we have therefore
\begin{equation}\label{eq-7}
(v-2) \left|G_{xB} \right|= (k-1)(k-2) \left|G_{xy}\right|,
\end{equation}
where
\[\left|G_{xB} \right|=n \frac{\left|G_{xy}
\right|}{\left|z^{G_{xy}}\right|}\] for some $n \in \N$. This is
equivalent to
\[2(2^{2d-2}\pm 2^{d-2}-1)n=(k-1)(k-2)\left| z^{G_{xy}} \right|\]
with
\[\left| z^{G_{xy}} \right| =\left\{\begin{array}{ll}
    2(2^{d-1} \mp 1)(2^{d-2} \pm 1),\;\mbox{or}\\
    2^{2d-2}.\\
\end{array} \right.\]
Clearly, $2^{2d-2}\pm 2^{d-2}-1 \equiv 1$ (mod $2$) and
$(k-1)(k-2) \equiv 0$ (mod $2$). As $(2^{2d-2}\pm 2^{d-2}-1,
2^{d-1} \mp 1) =(2^{d-2}, 2^{d-1} \mp 1)=1$ and $(2^{2d-2}\pm
2^{d-2}-1, \linebreak 2^{d-2} \pm 1)=( 2, 2^{d-2} \pm 1)=1$, it
follows that $\left| z^{G_{xy}} \right|$ always divides $n$. Thus
$\left| G_{xy} \right| \bigm| \left| G_{xB} \right|$, and
equation~(\ref{eq-7}) yields
\[v-2 \bigm| (k-1)(k-2).\]
But, on the other hand, we have $v-2 \geq (k-1)(k-2)$ by
Proposition~\ref{Cam}~(b), and it is immediately seen that $v$
cannot take the values where equality holds.

\bigskip
\emph{Cases} (7)-(8).
\medskip

For the existence of non-trivial flag-transitive Steiner
\mbox{$3$-designs}, we have in these cases only a small number of
possibilities for $k$ to check, which can easily be ruled out by
hand using Lemma~\ref{divprop}, Lemma~\ref{Comb_t=3}~(d), and
Corollary~\ref{Cameron_t=3}.

\pagebreak

\emph{Case} (9): $N=M_v$, $v=11,12,22,23,24$.
\medskip

Here $G$ is always \mbox{$3$-transitive}, and
thus~\cite[Thm.\,3]{Kant1985} yields the design described in part
(iv) of the Main Theorem. Obviously, 
flag-transitivity holds as the $3$-transitivity of $G$ implies that 
$G_x$ acts block-transitively on the derived Steiner 
\mbox{$2$-design} $\D_x$ for any $x \in X$.

\bigskip
\emph{Cases} (10)-(13).
\medskip

Again, the few possibilities for $k$ can easily be ruled out by
hand using Lemma~\ref{divprop}, Lemma~\ref{Comb_t=3}~(c) and (d),
and Corollary~\ref{Cameron_t=3}.

\medskip

\noindent This completes the proof of the Main Theorem.

\medskip


\subsection*{Acknowledgment}
I am grateful to C. Hering and W. M. Kantor for helpful
suggestions.

\bibliographystyle{amsplain}
\bibliography{Xbib3des}

\providecommand{\bysame}{\leavevmode\hbox to3em{\hrulefill}\thinspace}
\providecommand{\MR}{\relax\ifhmode\unskip\space\fi MR }
\providecommand{\MRhref}[2]{%
  \href{http://www.ams.org/mathscinet-getitem?mr=#1}{#2}
}
\providecommand{\href}[2]{#2}
\begin{thebibliography}{10}

\bibitem{Asch1987}
M.~Aschbacher, \emph{Chevalley groups of type {$G_2$} as the group of a
  trilinear form}, J. Algebra \textbf{109} (1987), 193--259.

\bibitem{Baer1946}
R.~Baer, \emph{Polarities in finite projective planes}, Bull. Amer. Math. Soc.
  \textbf{52} (1946), 77--93.

\bibitem{Beuk1981}
F.~Beukers, \emph{On the generalized {Ramanujan-Nagell} equation, {I}}, Acta
  Arith. \textbf{38} (1981), 389--410.

\bibitem{Block1965}
R.~E. Block, \emph{Transitive groups of collineations on certain designs},
  Pacific J. Math. {\bf 15} (1965), 13--18.

\bibitem{Buek1988}
F.~Buekenhout, A.~Delandtsheer, and J.~Doyen, \emph{Finite linear spaces with
  flag-transitive groups}, J. Combin. Theory, Series A \textbf{49} (1988),
  268--293.

\bibitem{Buek1990}
F.~Buekenhout, A.~Delandtsheer, J.~Doyen, P.~B. Kleidman, M.~W. Liebeck, and
  J.~Saxl, \emph{Linear spaces with flag-transitive automorphism groups}, Geom.
  Dedicata \textbf{36} (1990), 89--94.

\bibitem{Cam1976}
P.~J. Cameron, \emph{Parallelisms of {Complete Designs}}, London Math. Soc.
  Lecture Note Series {\bf 23}, Cambridge Univ. Press, Cambridge, 1976.

\bibitem{CaKa1979}
P.~J. Cameron and W.~M. Kantor, \emph{$2$-transitive and antiflag transitive
  collineation groups of finite projective and polar spaces}, J. Algebra
  \textbf{60} (1979), 384--422.

\bibitem{Atlas1985}
J.~H. Conway, R.~T. Curtis, S.~P. Norton, R.~A. Parker, and R.~A. Wilson,
  \emph{Atlas of {Finite} {Groups}}, Clarendon Press, Oxford, 1985.

\bibitem{CSK1976}
C.~W. Curtis, G.~M. Seitz, and W.~M. Kantor, \emph{The $2$-transitive
  permutation representations of the finite {Chevalley} groups}, Trans. Amer.
  Math. Soc. \textbf{218} (1976), 1--59.

\bibitem{Del1992}
A.~Delandtsheer, \emph{Finite (line, plane)-flag-transitive planar spaces},
  Geom. Dedicata \textbf{41} (1992), 145--153.

\bibitem{Del1995}
\bysame, \emph{Dimensional linear spaces}, in: Handbook of Incidence Geometry,
  ed. by F. Buekenhout, {Elsevier} {Science}, {Amsterdam}, 1995, 193-294.

\bibitem{Del2001}
\bysame, \emph{Finite flag-transitive linear spaces with alternating socle},
  in: Algebraic Combinatorics and Applications, Proc. Euroconf.
  (G\"{o}ßweinstein 1999), ed. by A. Betten et al., {Springer}, {Berlin}, 2001,
  79-88.

\bibitem{Deletal1986}
A.~Delandtsheer, J.~Doyen, J.~Siemons, and C.~Tamburini, \emph{Doubly
  homogeneous \mbox{$2$-$(v,k,1)$ designs}}, J. Combin. Theory, Series A
  \textbf{43} (1986), 140--145.

\bibitem{Demb1968}
P.~Dembowski, \emph{Finite {Geometries}}, Springer, Berlin, Heidelberg, New
  York, 1968; Reprint: Springer, 1997.

\bibitem{Dick1901}
L.~E. Dickson, \emph{Linear {Groups} with an {Exposition} of the {Galois Field
  Theory}}, Teubner, Leipzig, 1901; Reprint: Dover Publications, New York,
  1958.

\bibitem{Gor1982}
D.~Gorenstein, \emph{{Finite Simple Groups. An Introduction} to {Their
  Classification}}, Plenum Press, New York, London, 1982.

\bibitem{Her1974}
C.~Hering, \emph{Transitive linear groups and linear groups which contain
  irreducible subgroups of prime order}, Geom. Dedicata \textbf{2} (1974),
  425--460.

\bibitem{Her1985}
\bysame, \emph{Transitive linear groups and linear groups which contain
  irreducible subgroups of prime order, {II}}, J. Algebra \textbf{93} (1985),
  151--164.

\bibitem{Hu2001}
M.~Huber, \emph{Classification of flag-transitive {Steiner} quadruple systems},
  J. Combin. Theory, Series A \textbf{94} (2001), 180--190.

\bibitem{Hu2001diss}
\bysame, \emph{Klassifikationen fahnentransitiver {Steiner} {Designs}},
  Dissertation, Univ. T\"{u}bingen, T\"{u}bingen, 2001, (URL:
  http://w210.ub.uni-tuebingen.de/dbt/volltexte/ \linebreak 2001/275).

\bibitem{Hup1957}
B.~Huppert, \emph{Zweifach transitive, aufl\"{o}sbare {Permutationsgruppen}},
  Math. Z. \textbf{68} (1957), 126--150.

\bibitem{HupI1967}
\bysame, \emph{{Endliche} {Gruppen} {I}}, Springer, Berlin, Heidelberg, New
  York, 1967.

\bibitem{Kant1985}
W.~M. Kantor, \emph{Homogeneous designs and geometric lattices}, J. Combin.
  Theory, Series A \textbf{38} (1985), 66--74.

\bibitem{Kant1993}
\bysame, \emph{$2$-transitive and flag-transitive designs}, in: Coding Theory,
  Design Theory, Group Theory, Proc. Marshall Hall Conf., ed. by D. Jungnickel
  et al., {J. Wiley}, {New York}, 1993, 13-30.

\bibitem{KlLi1990}
P.~B. Kleidman and M.~W. Liebeck, \emph{The {Subgroup Structure} of the {Finite
  Classical Groups}}, London Math. Soc. Lecture Note Series {\bf 129},
  Cambridge Univ. Press, Cambridge, 1990.

\bibitem{Lieb1998}
M.~W. Liebeck, \emph{The classification of finite linear spaces with
  flag-transitive automorphism groups of affine type}, J. Combin. Theory,
  Series A \textbf{84} (1998), 196--235.

\bibitem{Luene1965}
H.~L\"{u}neburg, \emph{Fahnenhomogene {Quadrupelsysteme}}, Math. Z. \textbf{89}
  (1965), 82--90.

\bibitem{Mail1895}
E.~Maillet, \emph{Sur les isomorphes holo\'{e}driques et transitifs des groupes
  sym\'{e}triques ou altern\'{e}s}, J. Math. Pures Appl. (5) \textbf{1} (1895),
  5--34.

\bibitem{Saxl2002}
J.~Saxl, \emph{On finite linear spaces with almost simple flag-transitive
  automorphism groups}, J. Combin. Theory, Series A \textbf{100} (2002),
  322--348.

\bibitem{Suz1962}
M.~Suzuki, \emph{On a class of doubly transitive groups}, Ann. Math.
  \textbf{75} (1962), 105--145.

\bibitem{Ward1966}
H.~N. Ward, \emph{On {Ree's} series of simple groups}, Trans. Amer. Math. Soc.
  \textbf{121} (1966), 62--89.

\bibitem{Zsig1892}
K.~Zsigmondy, \emph{Zur {Theorie} der {Potenzreste}}, Monatsh. f\"{u}r Math. u.
  Phys. \textbf{3} (1892), 265--284.

\end{thebibliography}
\end{document}